\documentclass[11pt,reqno,twoside]{article}

\usepackage[hang]{footmisc}
\usepackage{lipsum}

\usepackage{fixltx2e} % Makes \( \) equation style robust, among other
\usepackage{cmap} % Load before fontenc 

\usepackage[T1]{fontenc}
\usepackage[utf8]{inputenc}
\usepackage{graphicx}
\usepackage{placeins}
\usepackage{enumerate}
\usepackage{algcompatible}
\usepackage{subcaption}
\usepackage{verbatim}
\newcommand{\comments}[1]{}

\usepackage{soul}

\usepackage{setspace}

\let\counterwithin\relax  %DSA: i had to include this to be able to compile
\usepackage{lmodern} % Improved version of computer modern
\usepackage[scale=0.88]{tgheros} % Helvetica clone for sans serif font

\usepackage{bm} % boldmath must be called after the package

\usepackage{bbold}

\usepackage{amsmath,amsbsy,amsgen,amscd,amsthm,amsfonts,amssymb} 

\usepackage[centering,top=1.0in,bottom=1.0in,left=0.75in,right=0.75in]{geometry}

\usepackage{titling}
\usepackage{musicography}
\graphicspath{ {./images/} }

\usepackage[sf,bf,compact]{titlesec}

\usepackage{booktabs,longtable,tabu} % Nice tables
\usepackage[font=small,margin=30pt,labelfont={sf,bf},labelsep={space}]{caption}

\usepackage[usenames,dvipsnames]{xcolor}

\definecolor{dark-gray}{gray}{0.3}
\definecolor{dkgray}{rgb}{.4,.4,.4}
\definecolor{dkblue}{rgb}{0,0,.5}
\definecolor{medblue}{rgb}{0,0,.75}
\definecolor{rust}{rgb}{0.5,0.1,0.1}

\usepackage{url}
\usepackage[colorlinks=true]{hyperref}
\hypersetup{linkcolor=dkblue}    
\hypersetup{citecolor=rust}      
\hypersetup{urlcolor=rust}     

\usepackage[final]{microtype} 

\newtheoremstyle{myThm} % name
    {\topsep}                    % Space above
    {\topsep}                    % Space below
    {\itshape}                   % Body font
    {}                           % Indent amount
    {\sffamily\bfseries}                   % Theorem head font
    {.}                          % Punctuation after theorem head
    {.5em}                       % Space after theorem head
    {}  % Theorem head spec (can be left empty, meaning ?€˜normal?€?)

\newtheoremstyle{myRem} % name
    {\topsep}                    % Space above
    {\topsep}                    % Space below
    {}                   % Body font
    {}                           % Indent amount
    {\sffamily}                   % Theorem head font
    {.}                          % Punctuation after theorem head
    {.5em}                       % Space after theorem head
    {}  % Theorem head spec (can be left empty, meaning ?€˜normal?€?)

\newtheoremstyle{myDef} % name
    {\topsep}                    % Space above
    {\topsep}                    % Space below
    {}                   % Body font
    {}                           % Indent amount
    {\sffamily\bfseries}                   % Theorem head font
    {.}                          % Punctuation after theorem head
    {.5em}                       % Space after theorem head
    {}  % Theorem head spec (can be left empty, meaning ?€˜normal?€?)

\theoremstyle{myThm}

\theoremstyle{myRem}

\theoremstyle{myDef}

\usepackage{fancyhdr}
\usepackage{algorithm}
\usepackage{algpseudocode}

\let\originalleft\left
\let\originalright\right
\renewcommand{\left}{\mathopen{}\mathclose\bgroup\originalleft}
\renewcommand{\right}{\aftergroup\egroup\originalright}

\usepackage{mathtools}
\mathtoolsset{centercolon}  % Makes := typeset correctly for definitions

\renewcommand{\phi}{\varphi}

\providecommand{\mathbbm}{\mathbb} % In case we don't load bbm

\newcommand{\R}{\mathbbm{R}}

\newcommand{\U}{\mathcal{U}}

\newcommand{\G}{\mathcal{G}}

\definecolor{mygreen}{rgb}{0.1,0.75,0.2}

\newcommand{\nc}{\normalcolor}

\newcommand{\M}{\mathcal{M}}

\newcommand{\Nc}{\mathcal{N}}

\usepackage[]{algorithm}

\usepackage{graphicx}

\usepackage{authblk}
\usepackage[square,numbers]{natbib}
\usepackage{chngcntr}
\usepackage{mathrsfs}
\counterwithin{table}{section}
\counterwithin{algorithm}{section}

\title{Graph-based  Prior and Forward Models for \\ Inverse Problems on Manifolds with Boundaries} % CORRECT TITLE
\author{John Harlim $^{1}$,  Shixiao Jiang $^{2},$ Hwanwoo Kim $^{3},$ and Daniel Sanz-Alonso $^{3}$}

\date{}

\vspace{.25in}

\makeatletter\@addtoreset{section}{part}\makeatother%
\numberwithin{equation}{section}

\newcommand{\upperRomannumeral}[1]{\uppercase\expandafter{\romannumeral#1}}

\begin{document}
\maketitle %  LEAVE HERE
% The command above causes the title to be displayed.

%>>>>> DELETE ALL CONTENT UNTIL "\end{document}"
% This is the body of your document.

\abstract{This paper develops manifold learning techniques for the numerical solution of PDE-constrained Bayesian inverse problems on manifolds with boundaries. We introduce graphical Mat\'ern-type Gaussian field priors that enable flexible modeling near the boundaries, representing boundary values by superposition of harmonic functions with appropriate Dirichlet boundary conditions. We also investigate the graph-based approximation of forward models from PDE parameters to observed quantities. In the construction of graph-based prior and forward models, we leverage the ghost point diffusion map algorithm to approximate second-order elliptic operators with classical boundary conditions. Numerical results validate our graph-based approach and demonstrate the need to design prior covariance models that account for boundary conditions. 
\let\thefootnote\relax\footnote{$^{1}$ Department of Mathematics, Department of Meteorology and Atmospheric Science, Institute for Computational and\\ \indent Data Sciences, The Pennsylvania State University, University Park, PA 16802, USA. }
\let\thefootnote\relax\footnote{$^2$ Institute of Mathematical Sciences, ShanghaiTech University, Shanghai 201210, China. }
\let\thefootnote\relax\footnote{$^3$ Department of Statistics, University of Chicago, Chicago, IL 60637, USA}

%\author{John Harlim, Shixiao Jiang, Hwanwoo Kim, Daniel Sanz-Alonso
%\\
%{\color{cyan}Here are our affiliations that you can arrange accordingly.}\\
%Shixiao Willing Jiang \\
%Institute of Mathematical Sciences, ShanghaiTech University, Shanghai
%201210, China\\ 
%\texttt{jiangshx@shanghaitech.edu.cn}  \vspace{5pt} \\
%John Harlim \\
%Department of Mathematics, Department of Meteorology and Atmospheric
%Science, \\
%Institute for Computational and Data Sciences \\
%The Pennsylvania State University, University Park, PA 16802, USA\\
%\texttt{jharlim@psu.edu}
%}

\section{Introduction}\label{sec:introduction}
This paper develops manifold learning techniques to address two tasks in the numerical solution of PDE-constrained Bayesian inverse problems on manifolds with boundaries: (1) the design and approximation of Gaussian field priors for spatially-distributed PDE parameters; and (2) the approximation of forward maps from PDE parameters to PDE solutions. We introduce graph-based approximations of prior and forward models and numerically show their effectiveness in two test problems. The first one concerns the recovery of the diffusion coefficient of an elliptic PDE from pointwise noisy measurements of the PDE solution; the second one concerns the recovery of the initial condition of a heat equation from noisy measurements of the PDE solution at some positive time. Both of these problems have been widely used as test cases for Bayesian inversion on manifolds and Euclidean domains \cite{AS10,franklin1970well,trillos2016bayesian,harlim2020kernel,trillos2020consistency,chada2020iterative,bigoni2020data}, but previous work has largely ignored the boundary effects that are the focus of this paper.

In the Bayesian approach to inverse problems \cite{tarantola2005inverse,kaipio2006statistical,AS10}, overviewed in Section \ref{sec: background} below, inference on the PDE parameters is performed using a posterior distribution obtained by conditioning a user-chosen prior distribution to  observed data. When the parameter of interest is a function, employing an adequate choice of prior is  crucial: the prior  determines the support of the posterior, and hence the space of parameters that can be recovered given sufficient data. In this work, we will focus on Mat\'ern-type Gaussian field priors  \cite{matern2013spatial,stein2012interpolation} that have been widely adopted in inverse problems and statistics due to their flexibility and  computational efficiency, e.g. \cite{lindgren2011explicit,bolin2019rational,bolin2020numerical, bolin2018weak,roininen2019hyperpriors,wiens2020modeling,bolin2014spatial}. We consider Mat\'ern-type models on manifolds \cite{lindgren2011explicit} and investigate their discretization using graphs, and their generalization to enhance flexibility near the boundaries. To our knowledge, all previous work on graph-based Mat\'ern models  \cite{bertozzi2018uncertainty,sanz2020spde,trillos2020consistency,sanz2020unlabeled}  disregarded boundary effects, despite  their importance in the discretization of Gaussian field priors in Euclidean settings \cite{khristenko2019analysis,stadler18,calvetti2006aristotelian}. Our numerical results will confirm that accounting for boundary conditions is important in the design of priors and in their graph-based approximation. While our emphasis is on PDE-constrained inverse problems, we expect that the boundary-aware graph-based Mat\'ern priors that we introduce will also find applications in graph-based machine learning \cite{sanz2020unlabeled}.

Approximation of the forward map requires solving PDEs on manifolds with a range of PDE parameters. There are numerous techniques to solve PDEs on manifolds, including the finite element method \cite{dziuk2013finite}, the level set method \cite{bertalmio2001variational}, closest point representation \cite{ruuth2008simple},  and the mesh-free radial basis function method \cite{fuselier2013high}. In contrast to all of these methods, the kernel approach we will consider avoids the need to have some parameterization on the manifolds; we refer to \cite{jiang2020ghost} for a detailed discussion of the advantages and disadvantages of these PDE solvers. In this work, we introduce graph-based forward map approximations for elliptic and heat inversion on manifolds with boundaries. Our graph-based approach leverages the ghost point diffusion map (GPDM) algorithm \cite{jiang2020ghost} to approximate second-order elliptic operators with classical boundary conditions, which are used to define these forward models and also the Mat\'ern prior covariance. The main idea behind the GPDM algorithm is to extend the domain of the PDE through a set of artificially constructed \emph{ghost points}. By extending the underlying manifold where the PDE was defined, one can treat the boundary of the original manifold as the interior of the extended manifold. The implication is that the kernel-based approximations that underpin the design of the GPMD method remain to be valid on the boundary points. The GPDM algorithm has been generalized to solve time-dependent advection-diffusion equations on manifolds, and the nonlinear viscous Burger's equation on a one-dimensional curve \cite{yan2021ghost}.

 We close this introduction with an outline of the paper and a summary of our contributions. 
\begin{itemize}
    \item In Section \ref{sec: background}, we overview the Bayesian approach to inverse problems and provide a high-level summary of the proposed procedure. We also review the Bayesian inversion approach in \cite{harlim2020kernel} for the elliptic inversion on a closed (i.e. boundary-free and compact) manifold. Lastly, we lay out the forward map discretization for the heat inverse problem.
    \item In Section \ref{sec:BEIPmb}, we review the GPDM algorithm and introduce our boundary-aware Mat\'ern-type priors and forward map approximations for elliptic and heat inversion.
    \item In Section \ref{sec:numerics}, we showcase implementations of the proposed methodology for elliptic and heat inversion on one and two-dimensional manifolds with boundaries. In addition, we demonstrate the enhanced flexibility of our proposed prior relative to existing graph-based approaches that ignore boundary effects.
    \item In Section \ref{sec:conclusions}, we provide conclusions and some open directions for future work.
\end{itemize}

\section{Background}\label{sec: background}
In this section, we first give the necessary background on the Bayesian formulation of inverse problems \cite{tarantola2005inverse,kaipio2006statistical,sanzstuarttaeb} and the function space perspective \cite{AS10, dashti2013bayesian}. Next, we formulate the problem of interest, namely, PDE-constrained Bayesian inversion on manifolds. Lastly, we review the graph-based approach in \cite{harlim2020kernel} which is only applicable for closed manifolds. 

\subsection{Bayesian Formulation of Inverse Problems}\label{subsec: bip}
%In this section we overview the Bayesian approach to inverse problems and introduce some notation. 
Suppose we have a parameter of interest $\theta \in \Theta$, observed data $y \in \mathbb{R}^M$ and a forward model $\mathcal{G},$ satisfying
\begin{equation}\label{eq:datamodel}
y = \mathcal{G}(\theta) + \eta,
\end{equation}
where $\eta \in \mathbb{R}^M$ is an observation noise. Our goal is to study the inverse problem of recovering $\theta$ from the observed data $y$. In the Bayesian framework, one endows probabilistic structure to any unknown quantity. Here this involves specifying a \emph{prior} distribution $\mu$ for the parameter $\theta$ and a distribution $\rho$ for the noise $\eta.$ To facilitate our presentation, we assume throughout that $\eta \sim \rho = \Nc( 0, \Gamma)$, where $\Gamma \in \mathbb{R}^{M \times M}$ is a given positive-definite covariance matrix. We also assume that $\theta$ and $\eta$ are independent random variables. The solution of the inverse problem under the Bayesian framework is then the \emph{posterior} distribution of $\theta$ given $y$, denoted $\mu^y$. The posterior distribution allows to not only construct point estimates of the parameter of interest but also quantify the uncertainty in the parameter reconstruction. In practice, this often requires to resort to sampling techniques such as Markov chain Monte Carlo (MCMC). 

 If the parameter space is finite-dimensional,  i.e.  $\Theta \subseteq \R^N,$ and the prior $\mu$ has Lebesgue density $\pi,$ then the posterior $\mu^y$ has Lebesgue-density $\pi^y$ given by
\begin{equation}\label{subsec:BIP:finiteconditional}
\pi^y(\theta) = \frac{1}{Z}\rho \bigl(y-\mathcal{G}(\theta)\bigr)\pi(\theta), \quad \quad Z := \int \rho \bigl(y-\mathcal{G}(\theta)\bigr)\pi(\theta) \, d\theta,
\end{equation}
where $\rho \bigl(y-\mathcal{G}(u)\bigr)$ is the data \emph{likelihood}. Under suitable assumptions \cite{AS10}, the characterization \eqref{subsec:BIP:finiteconditional} of the posterior can be generalized to infinite dimensional parameter space $\Theta$  by writing the posterior as a change of measure with respect to the prior:
\begin{align}\label{subsec:BIP:RND}
\begin{split}
\frac{d\mu^y}{d\mu}(\theta) = \frac{1}{Z} \, \rho \bigl(y-\mathcal{G}(\theta)\bigr) 
\propto \exp \bigl(-\Phi(\theta;y)\bigr),
\end{split}
\end{align}
where we have defined $\Phi(\theta;y) := \frac{1}{2} \|y-\mathcal{G}(\theta)\|_{\Gamma}^2$ and we set $\|\cdot \|_{\Gamma} := \|\Gamma^{-\frac{1}{2}}\cdot\|$.  
%In this paper we will employ function space MCMC \cite{beskos2008mcmc} algorithms with desirable scaling properties under discretization refinement for PDE-constrained inverse problems. 

\subsection{General Setting and Approach}\label{subsec: setup}
   We are interested in the inverse problem of recovering a parameter function $\theta \in \Theta$ of a PDE defined on a manifold $\mathcal{M} \subset \R^D$ from noisy measurements of the PDE solution $u\in \mathcal{U}$ at given locations  $\{\tilde x_m\}_{m=1}^M  \subset \mathcal{M}.$  Here both $\Theta$ and $\U$ are suitable function spaces.  The data model \eqref{eq:datamodel} is therefore given by 
   \begin{equation}\label{eq:generaldatamodel}
   y_m = u(\tilde x_m ) + \eta_m,  \quad m = 1, \ldots, M, \quad \quad \eta= \{\eta_m\}_{m=1}^M \sim \mathcal{N}(0, \Gamma),
   \end{equation}
with the \emph{forward model} $\mathcal{G}: \theta \mapsto \bigl(u(\tilde x_1), \ldots, u(\tilde x_M) \bigr)$  defined as the composition of a \emph{forward map} $\mathcal{F}: \theta \mapsto u$ from PDE input to PDE solution with an \emph{observation map}   $\mathcal{O} : u \mapsto  \bigl(u(\tilde x_1), \ldots, u(\tilde x_M) \bigr)$ from PDE solution to observed quantities. We will discretize the forward map  using a point cloud in $\mathcal{M}$, which is denoted by $\{x_n\}_{n=1}^N \supseteq \{\tilde x_m\}_{m=1}^M,$ where $M \le N.$ Our approach for the computational implementation of the inverse problem is then summarized in the following four steps:

\begin{enumerate}[(i)]
    \item {\bf Prior Specification}: Specify a prior distribution for the infinite dimensional parameter $\theta \in \Theta$.
    \item {\bf Prior and Forward Model Discretization}:
    \begin{itemize}
    \item    Discretize using graph-based techniques the prior distribution $\mu$ to obtain $\mu_N$, a prior distribution over $\theta_N := \bigl(\theta(x_1), \cdots, \theta(x_N)\bigr) \in \R^N$. Note that $\mu_N$ is a prior distribution over discrete functions $\theta_N$ defined on the point cloud $\{x_n\}_{n=1}^N$. 
    \item  Discretize using graph-based techniques the given forward map $\mathcal{F}: \theta \to u$ to obtain $\mathcal{F}_N: \theta_N \mapsto u_N \coloneqq \mathcal{F}_N(\theta_N) \in \mathbb{R}^N$, where $u_N$ is an approximation to the solution of the PDE evaluated along the point cloud, i.e., $u_N \approx \bigl(u(x_1), \ldots, u(x_N)\bigr)$. Furthermore, we define $\G_N(\theta_N) := \bigl( u_N(\tilde x_1), \ldots, u_N(\tilde x_M) \bigr) \in\mathbb{R}^M$ where $u_N(\tilde{x}_i)$ denotes the component of the vector $u_N\in \mathbb{R}^N$ that corresponds to $\tilde{x}_i\in \{x_n\}_{n=1}^N$.
\end{itemize}   
    \item  {\bf Sampling}: Use an MCMC algorithm to obtain samples from the posterior distribution over $\theta_N$, given by
    \begin{equation}\label{subsec:BIP:disRND}
\frac{d\mu_N^y}{d\mu_N}(\theta_N) \propto \exp(-\Phi_N(\theta_N;y)), ~\text{where}~ \Phi_N(\theta_N;y) \coloneqq \frac{1}{2} \|y-\mathcal{G}_N(\theta_N) \|^2_{\Gamma}.
\end{equation}
    \item {\bf Interpolation}: If desired, extend the samples to functions on $\mathcal{M}$ with an interpolation map.
\end{enumerate}
%To simplify the exposition, we will focus on the fully-supervised setting where $m = n,$ that is, observations are available along the entire point cloud $\{x_i\}_{i=1}^n$ and $\F_n = \G_n.$  We remark, however, that the methodology we propose is applicable in the semi-supervised setting where $m \ll n.$ 
This four-step approach was introduced in \cite{trillos2020consistency,garcia2018continuum} and previous works on the discretization process include \cite{harlim2020kernel}. The primary focus of this paper is to contribute to the specification and discretization steps for PDE-constrained inverse problems on manifolds with boundaries. We propose (1) priors that accommodate for prescribed boundary conditions, (2) graph-based discretizations of such priors; and (3) graph-based discretizations of the forward PDEs supplemented with boundary conditions. These procedures will be introduced in Section \ref{sec:BEIPmb}. In the remainder of this section, we address the sampling and interpolation steps, which are based on existing methodologies that are applicable in wide generality. 

\indent  For the sampling step, we will use the graph pCN method \cite{bertozzi2018uncertainty, trillos2020consistency}, summarized in Algorithm \ref{alg1}.
\begin{figure}
\makebox[\linewidth]{%
  \begin{minipage}{\dimexpr 0.75\linewidth+5em}
\begin{algorithm}[H]
\FloatBarrier
\caption{Graph pCN \label{alg1}}
\STATE {\bf Input}: initial value $\theta_N^{(0)}$, number of samples $J$, parameter value $\zeta \in (0,1).$ \\
\STATE {\bf For} $j = 1, \ldots, J$ {\bf do}: 
\begin{enumerate}
\item Propose $\tilde\theta_N^{(j)} = (1-\zeta^2)^{\frac{1}{2}}\theta_N^{(j)} + \zeta \, \xi_N^{(j)}$, where $\xi_N^{(j)} \sim \mathcal{N}(0, V_N)$.
\item  Set $$a(\theta_N^{(j)}, \tilde\theta_N^{(j)}) \coloneqq \min\left\{1, \exp\left(\Phi_N(\theta_N^{(j)};y) - \Phi_N(\tilde\theta_N^{(j)};y) \right)\right\},$$
and let
$$
\theta_N^{(j+1)} = \begin{cases} \tilde\theta_N^{(j)} &\text{with probability}~a(\theta_N^{(j)}, \tilde\theta_N^{(j)}), \\
\theta_N^{(j)} &\text{with probability}~1-a(\theta_N^{(j)}, \tilde\theta_N^{(j)}) .
\end{cases}
$$
\end{enumerate}

\STATE {\bf Output}: Samples of discrete functions $\{\theta_N^{(j)}\}_{j=1}^J.$
\end{algorithm}
  \end{minipage}}
\end{figure}
\noindent The covariance matrix $V_N$ in the pCN proposal is tightly connected with the prior construction. In fact, this matrix will be precisely the covariance matrix of the discretized prior distribution. The parameter $\zeta$ controls the size of the proposed moves of the chain. For a large $\zeta$ value, we explore a wider region of the state space with higher number of rejections. The motivation for using the graph pCN method is that it shows robustness with respect to the level of discretization refinement. We refer to \cite{beskos2008mcmc, cotter2013, bertozzi2018uncertainty, trillos2020consistency} for the theoretical and empirical justification of this robustness.

\indent For the interpolation step, once we obtain finite-dimensional samples $\theta_N \in \R^N$ representing a parameter function evaluated on the point cloud $\{x_n\}_{n =1}^N$, we can extend the samples into functions on $\M$ using the $\mathsf{K}$-NN interpolation map defined by
$$
\theta(x) = \frac{1}{\mathsf{K}}\sum_{x_i \in N_\mathsf{K}(x)} \theta_N(x_i), ~x \in \mathcal{M}, ~x \notin \{x_n\}_{n=1}^N,
$$
where $N_\mathsf{K}(x)$ is the set of $\mathsf{K}$-nearest points in the point cloud $\{x_n\}_{n=1}^N$ to the point $x$. To find the $\mathsf{K}$-nearest points, one can use the Euclidean distance in $\R^D$ or the geodesic distance on $\M \subset \R^D,$ if available. We refer to \cite{trillos2020consistency} for more details and we note that other interpolation methods are possible.

\subsection{Elliptic Inverse Problems and Heat Inversion on Closed Manifolds} \label{ssec:boundaryfreeIP}
In this section, we overview existing procedures for the prior specification and the discretization of prior and forward models on \emph{closed} (i.e. compact and boundary-free) manifolds. We start by describing the two inverse problems used as running examples. Throughout this section, $\mathcal{M}$ will denote a $d$-dimensional smooth closed manifold isometrically embedded in $\R^D.$ 

\paragraph{Elliptic Inverse Problem}
Consider the elliptic PDE
 \begin{equation}\label{subsec:BEIBF:inverse}
\mathcal{L}u \coloneqq -\text{div}(\kappa \nabla u) = f, ~x \in \mathcal{M},
\end{equation}
where the divergence and gradient operators are defined with respect to the Riemannian metric inherited by $\mathcal{M}$ from $\R^D.$ %from the $d$-dimensional smooth boundary-free manifold $\mathcal{M}$ embedded in $\mathbb{R}^D$. 
The goal of the elliptic inverse problem is to recover the diffusion coefficient $\kappa$ given the right-hand side $f$ (assumed to be smooth) and noisy pointwise observation of the solution $u$ at $M$ spatial locations $\{\tilde x_m\}_{m=1}^M \subset \M$.
The data are therefore given by  $\{y_m\}_{m=1}^M = \{ u(\tilde x_m) + \eta_m\}_{m=1}^M.$ 
 We cast this problem into the general framework \eqref{eq:datamodel} defining the forward map $\mathcal{F}: \theta \mapsto u$, where $\theta := \log \kappa \in (-\infty, \infty)$. Previous kernel-based methods for elliptic inverse problems on manifolds \cite{harlim2020kernel} assumed $\mathcal{M}$ to be closed so that the elliptic PDE given in \eqref{subsec:BEIBF:inverse} is not supplemented with a boundary condition. Bayesian elliptic inverse problems on Euclidean domains are one of the standard model problems in uncertanty quantification  \cite{AS10,dashti2013bayesian,trillos2016bayesian}. We will refer to the operator $\mathcal{L}$  as the weighted Laplacian operator. 

%The goal is to recover the diffusion coefficient $\kappa$ of an elliptic PDE defined along some manifold $\mathcal{M}$. More specifically, the elliptic PDE is given by
%\begin{equation}\label{subsec:BEIBF:inverse}
%\mathcal{L}u \coloneqq -\text{div}(\kappa \nabla u) = f, ~x \in \mathcal{M},
%\end{equation}
%where $f$ is assumed to be given function and smooth, and the divergence and gradient operators are defined with respect to the Riemannian metric inherited from the $m$-dimensional smooth manifold $\mathcal{M}$ embedded in $\mathbb{R}^d$. In addition, \cite{harlim2020kernel} assumed $\mathcal{M}$ to be compact and hence, the elliptic PDE given in \eqref{subsec:BEIBF:inverse} does not have any boundary condition. Corresponding to the PDE given by \eqref{subsec:BEIBF:inverse}, the forward map of interest is a mapping from the diffusion coefficient $\kappa$ to the solution $u$. In fact, \cite{harlim2020kernel} considered a forward map $\mathcal{F}: \theta \mapsto u$, where $\theta = \log \kappa \in (-\infty, \infty)$.

\paragraph{Heat Inversion} Consider the heat equation 
\begin{align}\label{subsubsec:BIheatEqn}
\begin{split}
\begin{cases}
u_t  = -\Delta_{\mathcal{M}} u , \quad \quad \quad \quad &(x, t) \in \mathcal{M}\times[0,\infty),  \\
u(x,0) = \theta(x),  & x \in \mathcal{M},
\end{cases}
\end{split}
\end{align}
where $\Delta_{\mathcal{M}} \coloneqq -\text{div}( \nabla \cdot)$ is the Laplace-Beltrami operator on $\mathcal{M}$. The goal of the inverse heat problem  is to recover the initial heat $\theta$ defined on $\mathcal{M}$ from noisy pointwise observation of the heat at time $t^* > 0$ along $M$ spatial locations $\{\tilde x_m\}_{m=1}^M \subset \M$. The data are therefore given by $\{y_m\}_{m=1}^M = \{u(\tilde x_m, t^*) + \eta_m\}_{m=1}^M$ where $u(\,\cdot\,,t^*)$ is the solution of the heat equation \eqref{subsubsec:BIheatEqn} at time $t^*$. The corresponding forward map is a heat equation solver, namely $\mathcal{F}: \theta \mapsto u$, where $u(x) := u(x, t^*)$  and $\theta$ is the initial condition. The function space formulation of Bayesian heat inversion was introduced in \cite{franklin1970well} on Euclidean domains, and graph-based formulations on closed manifolds were studied in \cite{trillos2020consistency}. Other than its natural physical interpretation, this inverse problem can also be viewed as a standard deconvolution problem arising in imaging applications.

\subsubsection{Prior Specification and Discretization}\label{subsec:priorspec}
Here we describe the specification of Mat\'ern Gaussian field priors for the log-diffusion coefficient and the initial heat in closed manifolds. We will also overview their graph-based discretization. Recall that the Mat\'ern model on a closed manifold $\M$ is defined as the Gaussian measure
\begin{equation}\label{subsec:BEIBF:prior}
\mu = \mathcal{N}(0, V), \quad  V = c (\tau I + \Delta_\mathcal{M})^{-s},
\end{equation}
where $\Delta_\mathcal{M}$ is the Laplace-Beltrami operator on $\mathcal{M}$ and $\tau > 0, s > \frac{d}{2}$ are two free parameters. The choice of normalizing constant 
$$
c = \frac{1}{\sum_{i=1}^\infty (\tau + \lambda_i)^{-s}},
$$
where $(\lambda_i)_{i=1}^\infty$ are the increasingly ordered eigenvalues of $\Delta_\mathcal{M},$ ensures that $\theta \sim \mu$ has unit marginal variance. 
Samples can be represented using the Karhunen-Lo\`eve (KL) expansion
\begin{equation}\label{subsec:BEIBF:KLprior}
\theta(x) = \sqrt{c}\sum_{i=1}^\infty(\tau + \lambda_i)^{-\frac{s}{2}}\zeta_i \phi_i(x), \quad x \in \M,
\end{equation}
where $(\zeta_i)_{i=1}^\infty$ are i.i.d. standard normal random variables and $(\phi_i)_{i=1}^\infty$ are eigenfunctions of $\Delta_\mathcal{M}$ with corresponding eigenvalues $(\lambda_i)_{i=1}^\infty$. The parameter $\tau$ represents an inverse length-scale and allows to discern the significant terms in the KL expansion  \eqref{subsec:BEIBF:KLprior}. The parameter $s$ characterizes the almost-sure regularity of the samples. The requirement $s > \frac{d}{2}$ is motivated by Weyl's law, which asserts that $\lambda_i \asymp i^{\frac{2}{d}}$, i.e., the asymptotic behavior of $\lambda_i$ is equivalent to that of $i^{\frac{2}{d}}$. Thus, $s>\frac{d}{2}$ ensures that $\mu$ is a well-defined Gaussian measure in $L^2(\mathcal{M}).$ 

To obtain the discretized prior distribution over the discrete functions defined along the point cloud $\{x_n\}_{n=1}^N \subset \mathcal{M}$, the paper \cite{harlim2020kernel} proposed to replace the Laplace-Betrami operator $\Delta_{\mathcal{M}}$ in \eqref{subsec:BEIBF:prior} via a graph Laplacian $\Delta_N$ which can be constructed using the point cloud $\{x_n\}_{n=1}^N$. To be more specific, the discretized prior distribution is given by
$$
\mu_N = \mathcal{N}(0, V_N), \quad  V_N = c_N(\tau I + \Delta_N)^{-s},
$$
where $\Delta_N \in \mathbb{R}^{N\times N}$ is a graph Laplacian  constructed using $\{x_n\}_{n=1}^N$ and $\tau > 0, \, s > \frac{d}{2}$ are two free parameters. The role of the parameters $\tau$ and $s$ is analogous to the infinite-dimensional case.  In practice, among different choices of graph Laplacian \cite{von2007tutorial},  the use of self-tuning graph Laplacian \cite{zelnik2005self}  was recommended in \cite{harlim2020kernel}. We refer to \cite{von2007tutorial} for a review of graph Laplacians and to \cite{sanz2020spde} for generalizations of this graph-based prior model and its connection to the Mat\'ern family \cite{stein2012interpolation}.
To obtain samples whose variance per node is one, the normalizing constant $c_N$ can be set to be 
$$
c_N = \frac{N}{\sum_{n=1}^N (\tau + \lambda_n^{(N)})^{-s}},
$$
where $(\lambda_n^{(N)})_{n=1}^N$ are the increasingly ordered eigenvalues of $\Delta_N$. Samples $\theta_N$ from this prior can be represented via the KL expansion
\begin{equation}\label{subsec:BEIBF:KLdprior}
\theta_N(x_i) = \sqrt{c_N}\sum_{n=1}^N(\tau + \lambda_n^{(N)})^{-\frac{s}{2}}\zeta_n \phi_n^{(N)}(x_i), \quad 1 \le i \le N,
\end{equation}
where $(\zeta_n)_{n=1}^N$ are i.i.d. standard normal random variables and $(\phi_n^{(N)})_{n=1}^N$ are eigenvectors of $\Delta_N$ with corresponding eigenvalues $(\lambda_n^{(N)})_{n=1}^N$. 

\subsubsection{Forward Map Discretization}
\paragraph{Elliptic Inverse Problem}
For the discretization of the elliptic forward map, one can approximate the operator $\mathcal{L}$ in equation \eqref{subsec:BEIBF:inverse} by an integral operator \cite{harlim2020kernel}, which can be subsequently approximated using a Monte Carlo  sum. To be more specific, let 
$$
G_{\epsilon} \bigl(u(x)\bigr) \coloneqq  \epsilon^{-\frac{d}{2}}\int_{\mathcal{M}} h \left( \frac{|x-y|^2}{4\epsilon}\right)u(y)  \, d\mathcal{V} (y),
$$
where $h(s) = 2^{-d}\pi^{-\frac{d}{2}} \exp(-s)$ and $\mathcal{V}$ denotes the volume form inherited by $\mathcal{M}$ from the ambient space $\R^D$. For a smooth $u$, it was shown in \cite{coifman2006diffusion} that
\begin{equation}\label{subsec:BEIBF:coifman}
G_{\epsilon} \bigl(u(x)\bigr) = u(x) + \epsilon \bigl(\omega u(x) - \Delta_{\mathcal{M}} u(x)\bigr) + \mathcal{O}(\epsilon^2), ~x \in \mathcal{M},
\end{equation}
where $\omega$ is a function that depends only on the parametrization of the manifold $\mathcal{M}$. %Combining \eqref{subsec:BEIBF:coifman} with the identity
Note that 
$$
\mathcal{L}u = -\kappa \Delta_{\mathcal{M}}u - \nabla u \cdot \nabla \kappa = \sqrt{\kappa} \Bigl(u \Delta_{\mathcal{M}}(\sqrt{\kappa}) - \Delta_{\mathcal{M}}(u\sqrt{\kappa})\Bigr),
$$
and use \eqref{subsec:BEIBF:coifman} on $\sqrt{\kappa}$ and $u\sqrt{\kappa}$ to obtain, 
\begin{align*}
    uG_{\epsilon}\sqrt{\kappa} &= u\sqrt{\kappa} + \epsilon \bigl(\omega u\sqrt{\kappa} - u\Delta_{\mathcal{M}}\sqrt{\kappa}\bigr) + \mathcal{O}(\epsilon^2) \\
    G_{\epsilon}(u\sqrt{\kappa}) &= u\sqrt{\kappa} + \epsilon \bigl(\omega u\sqrt{\kappa} - \Delta_{\mathcal{M}}\bigl(u\sqrt{\kappa}\bigr)\bigr) + \mathcal{O}(\epsilon^2).
\end{align*}
This yields
$$
G_{\epsilon}(u\sqrt{\kappa}) - uG_{\epsilon}\sqrt{\kappa} = \epsilon\Bigl(u \Delta_{\mathcal{M}}(\sqrt{\kappa}) - \Delta_{\mathcal{M}}(u\sqrt{\kappa})\Bigr) + \mathcal{O}(\epsilon^2) = \frac{\epsilon}{\sqrt{\kappa}}\mathcal{L}u + \mathcal{O}(\epsilon^2).
$$
Motivated from this equation, we define 
\begin{align*}
    \mathcal{L}_\epsilon u(x) \coloneqq \frac{\sqrt{\kappa(x)}}{\epsilon}\bigl(G_{\epsilon}(u(x)\sqrt{\kappa(x)}) - u(x)G_{\epsilon}\bigl(\sqrt{\kappa(x)}\bigr) \bigr),
\end{align*}
which can be rewritten as the following integral operator
\begin{equation}\label{subsec:BEIBF:kernelanisotropic}
\mathcal{L}_\epsilon u(x) \coloneqq \frac{1}{\epsilon^{\frac{d}{2}+1}} \int_{\mathcal{M}} h\left( \frac{|x-y|^2}{4\epsilon}\right) \sqrt{\kappa(x)\kappa(y)} \bigl(u(y)-u(x)\bigr) \, d {\mathcal{V}} (y)
\end{equation}
satisying
$$
\mathcal{L}_{\epsilon}u(x) = \mathcal{L}u(x) + \mathcal{O}(\epsilon), ~x \in \mathcal{M}.
$$
%For more detailed derivations, we refer to \cite{coifman2006diffusion,harlim2020kernel, harlim2018data}. 

The kernel operator $\mathcal{L}_\epsilon$  % defined in \eqref{subsec:BEIBF:kernelanisotropic}
 can be approximated by Monte Carlo viewing the point cloud $\{x_n\}_{n=1}^N$ as manifold samples, giving 
\begin{equation}\label{subsec:BEIBF:discreteanisotropic}
\mathcal{L}_{\epsilon,N} u(x_i) \coloneqq \frac{1}{\epsilon^{\frac{d}{2}+1}} \left(\frac{1}{N} \sum_{j=1}^N h\left( \frac{|x_i-x_j|^2}{4\epsilon}\right) \sqrt{\kappa(x_i)\kappa(x_j)} \,q_\epsilon(x_j)^{-1} \bigl(u(x_j)-u(x_i)\bigr)\right), ~i = 1, \ldots, N,
\end{equation}
where 
$$
q_\epsilon(x_j) := \frac{1}{2^{d}\pi^{\frac{d}{2}}N\epsilon^{\frac{d}{2}}}\sum_{k=1}^N \exp\left(-\frac{|x_j-x_k|^2}{4\epsilon}\right)
$$
is a kernel estimator of the density of the point cloud. One can write the discretized weighted Laplacian operator in \eqref{subsec:BEIBF:discreteanisotropic} in a matrix form.
 Specifically, define a kernel matrix $H$ with entries  $H_{ij} := \exp\left(-\frac{|x_i-x_j|^2}{4\epsilon}\right)$ and define a vector $Q $ with entries $Q_i =\sum_{j=1}^N H_{ij}$. Set the matrix $W$ with entries $W_{ij} = \sqrt{\kappa(x_i)\kappa(x_j)} \,H_{ij}Q_j^{-1}$ and the diagonal matrix $D$ with diagonal entries $D_{ii}=\sum_{j=1}^N W_{ij}$. Then the discretized weighted Laplacian $\mathcal{L}_{\epsilon,N}$ can be written as
\begin{equation}\label{subsec:weightedLaplace}
\mathcal{L}_{\epsilon,N} = \frac{1}{\epsilon}(W-D).
\end{equation}

For practical implementation, the bandwith parameter $\epsilon$ can be empirically chosen so that it lies in the region where  
\begin{equation}\label{subsec:epscalibration}
\log\bigl(T(\epsilon)\bigr) = \log\left(\sum_{i,j = 1}^{N, K} \exp\left(-\frac{|x_i - x_j(i)|^2}{4\epsilon} \right)\right)
\end{equation}
is approximately linear \cite{harlim2020kernel}. Here $(x_j(i))_{j=1}^K$ are the $K$ closest points to the point $x_i$. In practice, it was further observed that the maximum slope of $\log\left(T(\epsilon)\right)$ often coincided with $\frac{d}{2}$, where $d$ is the dimension of the underlying manifold $\mathcal{M}$ (see \cite{berry2016variable}).

Using the above discretization, one can obtain a discretized forward map, $\mathcal{F}_{\epsilon, N}$, which maps $\theta_N = \bigl(\log(\kappa(x_1)), \ldots, \log(\kappa(x_N))\bigr)$ to an $N$-dimensional vector $u_N$, which represents an approximate solution to \eqref{subsec:BEIBF:inverse} restricted to the point cloud. In other words, $u_N$ is the minimal norm least-squares solution of  % the following problem
$$
\mathcal{L}_{\epsilon, N}u_N = f_N,
$$
where $f_N = \bigl(f(x_1), \cdots, f(x_N)\bigr)$. Therefore, the discretized forward map is given by 
$$
\mathcal{F}_{\epsilon, N}: \theta_N \mapsto u_N = \mathcal{L}_{\epsilon, N}^{-1} f_N,
$$
where $\mathcal{L}_{\epsilon, N}^{-1}$ denotes the pseudo-inverse. This allows us to write the discretized posterior distribution $\mu^y_N$ as a change of measure with respect to the discretized prior distribution $\mu_N$ in the following way:
$$
\frac{d\mu^y_N}{d\mu_N}(\theta_N) \propto \exp\left(-\frac{1}{2}|y-\mathcal{G}_{\epsilon, N}(\theta_N)|_{\Gamma}^2\right),$$
where $\mathcal{G}_{\epsilon, N}(\theta_N) = \bigl( u_N(\tilde x_1), \ldots, u_N(\tilde x_M) \bigr).$ 
The Lebesgue density of the posterior is given by $\pi_N^y (\theta_N) \propto \rho \bigl( y - \mathcal{G}_{\epsilon, N}(\theta_N) \bigr) \pi_N(\theta_N)$, where $\rho$ and $\pi_N$ are Gaussian densities $\Nc(0,\Gamma)$ and $\Nc(0,V_N)$ in $\R^M$ and $\R^N$, respectively. 
The graph pCN algorithm can then be used to obtain samples $\theta_N \sim \pi_N^y$ that can be extended into the underlying manifold using a $\mathsf{K}$-NN map. 
%\indent \red Meanwhile, a wide range of real-world applications specifies the boundary conditions, and therefore, it is necessary to extend their methods to the elliptic boundary value problems. With this goal in mind, notice the key underlying assumption on the asymptotic expansion expression in \eqref{subsec:BEIBF:coifman} is that the point of evaluation $x$ must be away from the boundary by at least $\epsilon^r$ distance, for some small $\epsilon > 0$ and $r \in (0, \frac{1}{2})$ \cite{jiang2020ghost}. In \cite{harlim2020kernel}, the authors assumed $\mathcal{M}$ to be a compact manifold, thereby requiring no treatment of boundary conditions. This facilitated them to use the above asymptotic expression to construct the approximate operator of the anisotropic diffusion operator. Naturally, for the Bayesian elliptic inverse problem on a manifold with a boundary condition, we would need to circumvent the difficulty of using the above asymptotic expression in constructing a kernel approximated version of the anisotropic diffusion operator. Further details on this issue are covered in Section \ref{sec:BEIPmb}. [I believe this can be deleted from here and explained in the first paragraph of Section 3] \nc

\paragraph{Heat Inversion}
To discretize the heat forward map, we replace the Laplace-Beltrami operator in \eqref{subsubsec:BIheatEqn} with a graph Laplacian $\Delta_N$ and solve
\begin{align*}
\begin{cases}
\frac{\partial}{\partial t}u_N = -\Delta_N u_N, \quad \quad    ~t \in [0, \infty), \\
u_N(0) = \theta_N,
\end{cases}
\end{align*}
where  $\theta_N \in \R^N$ represents the initial heat function restricted to the point cloud.
Then the solution $u_N := u_N(t^*)$ at time $t^*$ can be expressed as 
$$
u_N= \sum_{n=1}^N \langle \theta_N , \phi_n^{(N)}\rangle e^{-\lambda_n^{(N)} t^*} \phi_n^{(N)}, 
$$
where $(\lambda_n^{(N)}, \phi_n^{(N)})_{n=1}^N$ are the ordered eigenpairs of the graph Laplacian $\Delta_N.$  From this, we naturally obtain a discretization of the forward map given by
$$
\mathcal{F}_N: \theta_N \mapsto u_N = \sum_{n=1}^N \langle \theta_N, \phi_n^{(N)}\rangle e^{-\lambda_n^{(N)} t^*} \phi_n^{(N)}.
$$
We can then write the discretized posterior distribution $\mu_N^y$ as a change of measure with respect to the prior distribution $\mu_N$
$$
\frac{d\mu^y_N}{d\mu_N}(\theta_N) \propto \exp\left(-\frac{1}{2}|y-\mathcal{G}_N(\theta_N)|_{\Gamma}^2\right),
$$
where $\G_N(\theta_N) = \bigl( u_N(\tilde x_1), \ldots, u_N(\tilde x_M) \bigr).$ 
The Lebesgue density of the posterior is given by $\pi_N^y (\theta_N) \propto \rho \bigl( y - \mathcal{G}_N(\theta_N) \bigr) \pi_N(\theta_N)$,  and sampling and interpolation can be performed with the general methodology described previously.

\section{Bayesian Inverse Problems on Manifolds with Boundaries}\label{sec:BEIPmb}

In this section, we propose a novel methodology for Bayesian inversion on manifolds with boundaries, addressing the design of priors that can  reflect  the given boundary conditions, the graph-based discretization of these priors, and the graph-based approximation of PDE-constrained forward maps supplemented with boundary conditions. To be concrete, we will focus on the following elliptic and heat inverse problems on one and two-dimensional manifolds with Dirichlet boundary conditions. Throughout this section and the rest of this paper, $\mathcal{M}$ will denote a smooth compact $d$-dimensional manifold isometrically embedded in $\mathbb{R}^D$  with boundary $\partial \mathcal{M}.$ We denote $\mathcal{M}^o = \mathcal{M}\setminus \partial \mathcal{M}.$ 

\paragraph{Elliptic Inverse Problem}
Consider the following elliptic PDE with Dirichlet boundary conditions, 
\begin{equation}\label{eq:ellipticPDEdirichlet}
\begin{cases}
\mathcal{L}u \coloneqq -\text{div}(\kappa \nabla u) = f, \quad x \in \mathcal{M}^o,\\
u(x) = h(x), \quad x \in \partial\mathcal{M},\\
\end{cases}
\end{equation}
%where $\mathcal{M}$ is a $d$-dimensional manifold embedded in $\mathbb{R}^D$, $\partial \mathcal{M}$ is its boundary and $\mathcal{M}^o = \mathcal{M}\setminus \partial \mathcal{M}.$ 
where the divergence and gradient operators are defined with respect to the Riemannian metric inherited by the manifold $\mathcal{M}$ from $\mathbb{R}^D$. In our numerical examples, we will consider a one-dimensional semi-ellipse where  $\partial \mathcal{M}=\{x_0,x_N\}$ and a two-dimensional semi-torus where $\partial \mathcal{M}=\mathcal{B}_1\cup \mathcal{B}_2$ and $\mathcal{B}_i,$ $i =1,2$ are circles. We will then denote by $h_1$ and $h_2$ the Dirichlet boundary conditions at $\mathcal{B}_1$ and $\mathcal{B}_2.$
The goal of the elliptic inverse problem is to recover the diffusion coefficient $\kappa$ given the right-hand side $f$ and noisy pointwise observation of the solution $u$ at $M$ spatial locations $\{\tilde x_m\}_{m=1}^M \subset \mathcal{M}$. The data are given by $\{y_m\}_{m=1}^M = \{ u(\tilde x_m) + \eta_m\}_{m=1}^M$  with the forward map $\mathcal{F}: \theta \mapsto u$, where $\theta := \log \kappa \in (-\infty, \infty)$. 
% {\color{cyan} Is $\theta$ in the boundary condition in (3.1) related to $\log\kappa$? }

\paragraph{Heat Inversion} 
Consider the following one-dimensional heat equation with Dirichlet boundary conditions,
\begin{align}\label{eq:heatdirichlet}
\begin{cases}
u_t(x,t) = -\Delta_\mathcal{M} u (x,t),  \quad x \in \mathcal{M}^o, t > 0,   \\
u(x, t) = h(x), \quad ~x \in \partial\mathcal{M},\, t\ge 0, \\
u(x, 0) = \theta(x),  \quad x \in \mathcal{M},
\end{cases}
\end{align}
%Consider the following one-dimensional heat equation with Dirichlet boundary conditions,
%\begin{align}\label{subsec:BIPMBheat}
%\begin{cases}
%u_t(x,t) = \Delta_\mathcal{M}(x,t),  \quad x \in \mathcal{M}^o, t > 0,   \\
%u(x_1, t) = \theta(x_1), \quad t\ge 0, \\
%u(x_{n}, t) = \theta(x_{n}), \quad t \ge 0, \\
%u(x, 0) = \theta(x),  \quad x \in \mathcal{M},
%\end{cases}
%\end{align}
%where $\mathcal{M}, \mathcal{M}^o, \partial \mathcal{M}$ are defined as in the one-dimensional elliptic inverse problem. We also consider the two-dimensional heat equation with Dirichlet boundary conditions
%\begin{align}\label{subsec:BIPMB2Dheat}
%\begin{cases}
%u_t(x,t) = \Delta_\mathcal{M}(x,t),  \quad x \in \mathcal{M}^o, t > 0,   \\
%u(x, t) = h_1(x), ~x \in \mathcal{B}_1~t \ge 0, \\u(x, t) = h_2(x), ~x \in \mathcal{B}_2~t \ge 0, \\
%u(x, 0) = \theta(x), ~x \in \mathcal{M},
%\end{cases}
%\end{align}
%where the conditions and definitions for $\mathcal{B}_1, \mathcal{B}_2, h_1, h_2$, $\mathcal{M}$, $\partial \mathcal{M}$ are as in the two-dimensional elliptic inverse problem. 
Again the goal of the inverse heat problem is to recover the initial heat $\theta$ defined on $\mathcal{M}$ from noisy pointwise observation of the heat at time $t^* > 0$  along $M$ spatial locations $\{\tilde x_m\}_{m=1}^M \subset \M$. The data are therefore given by $\{y_m\}_{m=1}^M = \{u(\tilde x_m, t^*) + \eta_m\}_{m=1}^M$.  The corresponding forward map is a heat equation solver, namely $\mathcal{F}: \theta \mapsto u$, where $u(x) \equiv u(x, t^*)$ denotes the solution of \eqref{eq:heatdirichlet} at time $t^*$ and $\theta$ is the initial condition. For our numerical examples we will use a one-dimensional semi-ellipse and a two-dimensional semi-torus, and we will adopt the same notations as in the elliptic inverse problem. 
% {\color{cyan} Same comments and questions as above}

\subsection{Ghost Point Diffusion Maps for Dirichlet Boundary Conditions}\label{sec31}
In this section, we give a short discussion on the construction of ghost points and the ghost point diffusion map (GPDM) algorithm, which will be used in the discretizations of the prior and forward models. In this paper, we focus on 
a specific GPDM algorithm to approximate the weighted Laplacian operator, $\mathcal{L} =-\text{div}(\kappa \nabla \cdot)$ whose inputs are functions $u:\mathcal{M}\to \mathbb{R}$ that satisfy the Dirichlet boundary condition, $u(x)= h(x),$ for all $x\in\partial\mathcal{M}$. 
%Here,  $\mathcal{M}$ is a manifold with boundary, isometrically embedded in $\mathbb{R}^D$. 
The discussion will focus on the algorithmic aspect that will be used in the forward map discretization. For other types of (possibly non-symmetric) second-order diffusion operators and boundary conditions, and the convergence analysis, we refer to \cite{jiang2020ghost,yan2021ghost}. 

The key idea of the GPDM comes from the classical ghost point method \cite{leveque2007finite} for solving the PDE with Neumann boundary condition using the finite-difference method. Particularly, the ghost points are constructed to improve the convergence rate in approximating the normal derivative at the boundary points. In the context of the GPDM, the ghost points are employed to overcome the biases induced by the graph Laplacian discretization near the boundary. Numerically, solving PDEs with the ghost point method requires the following two steps: (1)  specification of ghost points; and (2) specification of function values at the ghost points. While these two steps are trivial when the geometry is Euclidean or known, they require nontrivial numerical algorithms and theoretical justification when the manifold is unknown in the sense that it can only be identified with finitely sampled point cloud data. 

The GPDM addresses step (1) above by augmenting the sampled point cloud data on the manifold with a set of ghost points specified on the outer normal collar of the boundary. Theoretically, the GPDM extends the embedded manifold $\mathcal{M} \stackrel{\iota}{\hookrightarrow} \mathbb{R}^D$ with its collar neighbor $\Delta\mathcal{M}$ of a sufficiently large radius, such that the extended manifold $\mathcal{M}\cup \Delta\mathcal{M} \stackrel{\tilde\iota}{\hookrightarrow} \mathbb{R}^D$ is isometrically embedded and does not change the geometry of $\mathcal{M}$, i.e., $\tilde{\iota}\vert_{\mathcal{M}}=\iota$ (see Lemma~3.5 in \cite{jiang2020ghost}). With this modification, the graph Laplacian construction in \eqref{subsec:BEIBF:discreteanisotropic} is a consistent pointwise estimator of $\mathcal{L}$ for all points in $\mathcal{M}$, even for the points that are very close to the boundary $\partial \mathcal{M}$ since they are sufficiently far away from the boundary of the extended manifold, $\partial (\mathcal{M}\cup \Delta \mathcal{M})$, as illustrated in Figure \ref{Ghost Point}. Since our goal is to construct a forward map on the manifold $\mathcal{M}$, we need to specify the additional unknowns (the function values at these ghost points as we noted in the step (2)) by adding more equations. Specifically, we will impose a set of linear extrapolation equations, whose solution specifies the function values on the ghost points through function values on the point cloud.

\begin{figure}
\centering
\includegraphics[width=.6\textwidth]{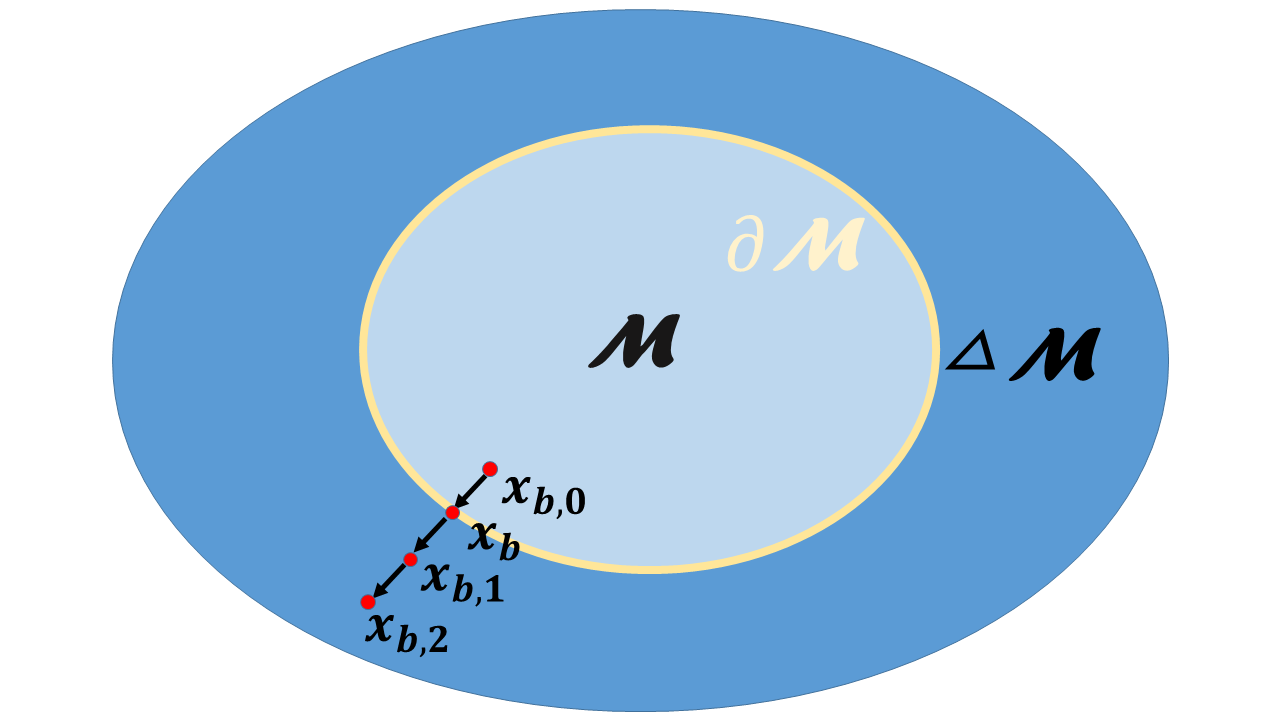}
\caption{Construction of ghost points.}
\label{Ghost Point}
\end{figure}

While there are various ways to realize these two steps above, in the following, we will provide a simple numerical procedure for well-ordered data that is used in the numerical examples in this paper. For randomly sampled data or a higher convergence rate method, we refer interested readers to \cite{jiang2020ghost}. We should also point out that the presentation below uses much simpler notations compared to those in \cite{jiang2020ghost,yan2021ghost} since we are only interested in the algorithmic aspect of the GPDM method. Given a point cloud $\{x_n\}_{n=1}^N\subset \mathcal{M}$, the GPDM estimator for $\mathcal{L}$ is constructed as follows:

\begin{enumerate}
\item Specification of the ghost points. At each boundary point, $x_b\in \partial \mathcal{M}$:  
\begin{description}
\item[(a)] Apply the secant line to approximate the normal vector at each boundary point. Following the notation in the illustration in Figure \ref{Ghost Point}, the unit normal vector is estimated via 
\begin{equation}
v_{b} = \frac{x_b - x_{b,0}}{\|x_b - x_{b,0}\|} .\notag
\end{equation}
\item[(b)] We specify $K$ ghost points along the normal vector at each boundary point, $x_b\in \partial \mathcal{M} \cap \{x_n\}_{n=1}^N$, 
$$
x_{b,k} \coloneqq x_b + h k v_{b}, \quad k = 1, \ldots, K,
$$
where $h=\|x_b - x_{b,0}\|$ and $x_{b,0}$ is the closest point to $x_b$ in $\mathcal{M}^o \cap \{x_n\}_{n=1}^N.$ 
\end{description}
\item Specification of the function values at the ghost points. For each boundary point $x_b$, we impose the following extrapolation condition,
\begin{equation}
\begin{aligned} u(x_{b,1}) - 2u(x_b)  +  u(x_{b,0}) &= 0,\\
u(x_{b,2})-2u(x_{b,1}) +u(x_b)&=0, \\ 
u(x_{b,k})-2u(x_{b,k-1})+ u(x_{b,k-2})&=0, \quad\quad k = 3,\ldots, K. \\
\end{aligned}  \label{Eqn:uvv_g2}
\end{equation}
These algebraic equations are  discrete analogs of matching the first-order derivatives along the estimated normal direction, $v_{b}$. 
\item Construction of the GPDM discrete estimator. Construct the graph Laplacian estimator as in \eqref{subsec:BEIBF:discreteanisotropic}
for the extended points, $\{x_n\}_{n=1}^N\cup\{x_{b,k}\}_{b,k=1}^{B,K}$. We should point out that we also use a set of algebraic equations identical to \eqref{Eqn:uvv_g2} to determine the function value of $\kappa$ at the ghost points $\{x_{b,k}\}_{b,k=1}^{B,K}$. Let $\bar{N}=N+BK$ and obtain the corresponding weighted Laplacian matrix, $\mathcal{L}_{\epsilon,\bar{N}}$, whose construction is analogous to \eqref{subsec:weightedLaplace}. The weighted Laplacian matrix can be rewritten as 
\begin{equation}
\mathcal{L}_{\epsilon,\bar{N}} = \begin{pmatrix} \mathcal{L}^{(1)}_{\epsilon,\bar{N}} & \mathcal{L}^{(2)}_{\epsilon,\bar{N}} \\
\mathcal{L}^{(3)}_{\epsilon,\bar{N}} & \mathcal{L}^{(4)}_{\epsilon,\bar{N}}
\end{pmatrix} \in \mathbb{R}^{\bar{N}\times\bar{N}}.
\end{equation}

We will construct the GPDM matrix based on the sub-matrices  $\mathcal{L}^{(1)}_{\epsilon,\bar{N}} \in \mathbb{R}^{N\times N}$, whose components are constructed based on the affinity between pairs of the point cloud points in $\{x_n\}_{n=1}^N\subset \mathcal{M}$, and $\mathcal{L}^{(2)}_{\epsilon,\bar{N}} \in \mathbb{R}^{N\times BK}$, whose components  represent  the affinity between an element of $\{x_n\}_{n=1}^N$ and a ghost point in $\{x_{b,k}\}_{b,k=1}^{B,K}$. Let $u_G = (u(x_{1,1}),\ldots,u(x_{B,K}))\in\mathbb{R}^{BK},$ $u_N = (u(x_1),\ldots,u(x_{b,0}) ,\ldots,u(x_N))\in\mathbb{R}^N$. Then we can write the solution of \eqref{Eqn:uvv_g2} in a compact form as $u_G= G u_N$ such that $\mathcal{L}^{(1)}_{\epsilon,\bar{N}}u_N  + \mathcal{L}^{(2)}_{\epsilon,\bar{N}} u_G = (\mathcal{L}^{(1)}_{\epsilon,\bar{N}}  + \mathcal{L}^{(2)}_{\epsilon,\bar{N}} G)u_N$. Based on this observation, we define the GPDM estimator without boundary condition as,
\begin{equation}
\tilde{\mathcal{L}}_{\epsilon,N}:= \mathcal{L}^{(1)}_{\epsilon,\bar{N}}  + \mathcal{L}^{(2)}_{\epsilon,\bar{N}} G \in \mathbb{R}^{N\times N}.
\end{equation}
\end{enumerate}

%Suppose we are given the following elliptic PDE with Dirichlet boundary conditions,
%\begin{equation}\label{subsec:EllipticInverseBoundary}
%\begin{cases}
%\mathcal{L}u \coloneqq -\text{div}(\kappa \nabla u) = f, ~x \in \mathcal{M}^o, \\
% u = h, ~x \in \partial{\mathcal{M}},
% \end{cases}
%\end{equation}
%where $\mathcal{M}^o, \partial{\mathcal{M}}$ respectively represent interior/boundary of the manifold $\mathcal{M},$ and $f, h$ are smooth functions on $\M$ and $\partial M,$ respectively.  
%{\color{cyan} We can erase the above equation and maybe just refer to (3.3) if you agree to merge (3.1) and (3.2) into (3.3).}
Denoting
\begin{equation}
\tilde{\mathcal{L}}_{\epsilon,N} = \begin{pmatrix} \tilde{\mathcal{L}}_{\epsilon,(N-B)\times(N-B)} &\tilde{\mathcal{L}}_{\epsilon,(N-B)\times B}  \\  
\tilde{\mathcal{L}}_{\epsilon,B\times(N-B)} &\tilde{\mathcal{L}}_{\epsilon,B\times B},
\end{pmatrix}
\end{equation}
and splitting $u_N:=(u_{N-B},u_B)$ into function values at the interior and boundary points, respectively, the GPDM discretizes the Dirichlet problem in \eqref{eq:ellipticPDEdirichlet} as follows,
\begin{eqnarray}\label{GPDMDirichlet}
\begin{aligned}
\tilde{\mathcal{L}}_{\epsilon,(N-B)\times(N-B)} {u}_{N-B} + \tilde{\mathcal{L}}_{\epsilon,(N-B)\times B} u_{B} &=& f_{N-B}, \\
u_B &=& h_B. 
\end{aligned}
\end{eqnarray}
Here, the components of the vector $f_{N-B} \in\mathbb{R}^{N-B}$ are the function values at the interior points, whereas the components of the vector $h_B\in \mathbb{R}^B$ are the function values at the boundary points.

\subsection{Prior Specification and Discretization}\label{subsec:BIEPHeatDirichlet}
In this section, we propose novel prior constructions to facilitate elliptic and heat Bayesian inversion on manifolds with boundaries. The priors we propose contain two terms. The first one is a Mat\'ern-type Gaussian field with Dirichlet homogeneous boundary conditions, whose primary role is to capture the uncertainty on the parameter of interest in the interior of the manifold. The second term accounts for prior uncertainty along the boundary of the manifold, and can be omitted if the boundary values of the parameter are known.  

For a one-dimensional manifold $\mathcal{M}$ with boundary $\partial \mathcal{M} = \{x_1, x_N\}$,  i.e.  a semi-ellipse in our numerical experiments, we propose using a prior defined as the law of 
\begin{equation}\label{subsubsec:1DPrior}
 \theta =    \underbrace{\frac{1}{\sqrt{\sum_{i=1}^\infty (\tau + \lambda_i)^{-s}}}\sum_{i=1}^\infty(\tau + \lambda_i)^{-\frac{s}{2}}\zeta_i \phi_i}_{\text{Interior term}} + \underbrace{\mu_1 \psi_1 + \mu_2\psi_2}_{\text{Boundary term}}.
\end{equation}
In the the interior term, $\tau, s,$ and $(\zeta_i)_{i=1}^\infty$ play the same role as in \eqref{subsec:BEIBF:KLprior}, but now $(\lambda_i, \phi_i)_{i=1}^\infty$ are the ordered eigenpairs of the Laplace-Beltrami operator $\Delta_{\mathcal{M}}$ equipped with homogeneous Dirichlet boundary conditions, so that $\phi_i(x_1) = \phi_i(x_N) = 0$ for all $i$.  Nonhomogeneous Dirichlet boundary conditions could also be considered. 
In the boundary term, $\mu_1, \mu_2 \sim \mathcal{N}(0,1)$ are independent of all other randomness and $\psi_1, \psi_2$ are harmonic functions defined by
\begin{align}\label{subsubsec:BIPMBnonhom}
\begin{cases}
\Delta_{\mathcal{M}} \psi_1 = 0,  \\
\psi_1(x_1) = 1, \\
\psi_1(x_N) = 0,
\end{cases}
\quad
\text{and}
\quad \quad
\begin{cases}
\Delta_{\mathcal{M}} \psi_2 = 0, \\
\psi_2(x_1) = 0, \\
\psi_2(x_N) = 1. 
\end{cases}
\end{align}
Thus the random coefficients $\mu_1$ and $\mu_2$ in \eqref{subsubsec:1DPrior}  respectively  represent the uncertainty on the parameter value of interest at the two boundary points $x_1$ and $x_N$. Overall, the prior defined by \eqref{subsubsec:1DPrior} combines the flexibility of the Mat\'ern model with homogeneous  Dirichlet boundary conditions for interior reconstruction with additional flexibility in the boundary reconstruction. A simulation study that illustrates the increased flexibility afforded by the boundary term and our prior construction is given in Section \ref{sec:numerics}. 

Following the same idea, for a two-dimensional manifold $\mathcal{M}$ with boundary  $\partial \mathcal{M} = \mathcal{B}_1 \cup \mathcal{B}_2$, where $\mathcal{B}_1$ and $\mathcal{B}_2$ are disjoint regular closed curves,  e.g. in our numerical example, $\mathcal{B}_1$ and $\mathcal{B}_2$ are boundary circles of a semi-torus, we define a prior as the law of
\begin{equation}\label{subsubsec:2DheatPrior}
 \theta  =   \underbrace{\frac{1}{\sqrt{\sum_{i=1}^\infty (\tau + \lambda_i)^{-s}}}\sum_{i=1}^\infty(\tau + \lambda_i)^{-\frac{s}{2}}\zeta_i \phi_i}_{\text{Interior term}} + \underbrace{\sum_{\ell=1}^L \mu_{1,\ell} \psi_{1,\ell} + \sum_{\ell=1}^L \mu_{2,\ell} \psi_{2,\ell}}_{\text{Boundary term}}.
\end{equation}
The interior term is defined analogously as in the one-dimensional case, using the spectrum of the Laplace-Beltrami operator with homogeneous Dirichlet boundary conditions. For the boundary term,
 $\mu_{1,\ell}, \mu_{2,\ell} \sim \mathcal{N}(0,1)$ are independent of all other randomness, and $\{\psi_{1,\ell}\}_{\ell=1}^L$ and $\{\psi_{2,\ell}\}_{\ell=1}^L$ are harmonic functions on $\mathcal{M}$ defined by 
\begin{align}\label{subsubsec:BIPMB2Dnonhom}
\begin{cases}
\Delta_{\mathcal{M}} \psi_{1,\ell} = 0,  \\
\psi_{1,\ell}(x) = \tilde\psi_{1,\ell}(x), \quad x \in \mathcal{B}_1, \\
\psi_{1,\ell}(x) = 0, \quad x \in \mathcal{B}_2,
\end{cases}
\quad
\text{and}
\quad \quad
\begin{cases}
\Delta_{\mathcal{M}} \psi_{1,\ell} = 0,  \\
\psi_{2,\ell}(x) = 0, \quad x \in \mathcal{B}_1, \\
\psi_{2,\ell}(x) = \tilde\psi_{2,\ell}(x), \quad x \in \mathcal{B}_2,
\end{cases}
\text{for} ~\ell = 1, \ldots, L.
\end{align}
where $ \tilde\psi_{1,\ell}(x)$ and $\tilde\psi_{2,\ell}(x)$ are ordered eigenfunctions of the Laplace Beltrami operator defined on the closed curves $\mathcal{B}_1$ and $\mathcal{B}_2,$ respectively. The number $L$ of basis-type functions for each boundary controls the flexibility of the prior along the boundary. Larger $L$ allows to recover more frequencies of the parameter of interest along the boundary, but at the expense of introducing additional model parameters.

To discretize the prior defined in \eqref{subsubsec:1DPrior}, we simply replace the role of the Laplace-Betrami operator $\Delta_{\mathcal{M}}$ with a graph Laplacian as we did in Subsection  \ref{subsec:priorspec}, but now taking care of boundary conditions. Specifically, for one-dimensional inverse problems, our discretized prior distribution would have samples of the form 
\begin{equation}\label{subsubsec:1DheatDPrior}
 \theta_N \sim   \frac{\sqrt{N}}{\sqrt{\sum_{n=1}^N (\tau + \lambda_n^{(N)})^{-s}}}\sum_{n=1}^N(\tau + \lambda_n^{(N)})^{-\frac{s}{2}}\zeta_n \phi_n^{(N)} + \mu_1 \psi_1^{(N)} + \mu_2\psi_2^{(N)},
\end{equation}
where $(\lambda_n^{(N)}, \phi_n^{(N)})$ are ordered eigenpairs of $\tilde \Delta_N$ defined as the $N\times N$ submatrix corresponding to the point cloud portion of the graph Laplacian constructed using both the point cloud $\{x_n\}_{n=1}^N$ and ghost points $\{x_{b,k}\}_{b,k=1}^{B,K}$. From now on, we will refer to $\tilde \Delta_N$ as a \emph{truncated graph Laplacian}. Notice that it provides an effective approximation to the Laplace Beltrami operator on $\M$ with homogeneous Dirichlet boundary conditions.  This construction avoids any potential idiosyncratic boundary behavior of eigenvectors when using a graph Laplacian constructed solely with point cloud data.  Indeed, in our empirical experience, eigenvectors of the self-tuned graph Laplacian tend to have spikes and oscillations near the boundary. Finally, $\psi_1^{(N)}$ and $\psi_2^{(N)}$ are the solutions of \eqref{subsubsec:BIPMBnonhom} evaluated along the point cloud. These solutions can be obtained by the GPDM as described in Section \ref{sec31}. 

For the two-dimensional graph-based prior discretization, we obtain a finite number of discretized functions for each boundary, denoted by $\{\psi_{1,\ell}^{(N)}\}_{\ell=1}^L$ and $\{\psi_{2,\ell}^{(N)}\}_{\ell=1}^L$, by solving \eqref{subsubsec:BIPMB2Dnonhom} along the point cloud. One can again use the GPDM to obtain $\{\psi_{1,\ell}^{(N)}\}_{\ell=1}^L$ and $\{\psi_{2,\ell}^{(N)}\}_{\ell=1}^L$ where the boundary basis functions $\{\tilde\psi_{1,\ell}\}_{\ell=1}^L$ and $\{\tilde\psi_{2,\ell}\}_{\ell=1}^L$ are discretized using a self-tuned graph Laplacian constructed with all the point cloud elements that lie on the boundaries, which are two disjoint closed curves in our numerical example. For instance, to obtain $\{\tilde\psi_{1,\ell}\}_{\ell=1}^L$, one can construct the self-tuned graph Laplacian solely using points in $\mathcal{B}_1 \cap \{x_n\}_{n=1}^N$, where $\{x_n\}_{n=1}^N \subset \mathcal{M}$ is the point cloud of the manifold. Then $\{\tilde\psi_{1,\ell}\}_{\ell=1}^L$ can be chosen to be the first $L$ eigenvectors of this self-tuned graph Laplacian. 
To summarize, in the two-dimensional case, samples from the proposed graph-based prior are defined by
\begin{equation}\label{subsubsec:2DPrior}
  \theta_N \sim   \frac{\sqrt{N}}{\sqrt{\sum_{n=1}^N (\tau + \lambda_n^{(N)})^{-s}}}\sum_{n=1}^N(\tau + \lambda_n^{(N)})^{-\frac{s}{2}}\zeta_n \phi_n^{(N)} +\sum_{\ell=1}^L \mu_{1,\ell} \psi_{1,\ell}^{(N)} + \sum_{\ell=1}^L\mu_{2,\ell} \psi_{2,\ell}^{(N)},
\end{equation}
where $(\lambda_n^{(N)}, \phi_n^{(N)})_{n=1}^N$ are ordered eigenpairs of a \emph{truncated graph Laplacian}.

%\blue 
%The use of the truncated graph Laplacian can be backed up by the following two aspects. Firstly, it avoids any potential idiosyncratic boundary behaviors of eigenvectors when using a graph Laplacian constructed solely with point cloud data. Under the Dirichlet boundary conditions, there is no established result of the convergence of graph Laplacian to the Laplace-Beltrami operator. Indeed, in our empirical experience, eigenvectors of the self-tuned graph Laplacian had peculiar behavior near the boundary of the underlying domain. Secondly, the truncated graph Laplacian, its eigenvectors tend to have boundary values close to zero. This second aspect is particularly useful as the first term in the discretized prior distribution are supposed to reflect the parameter of interest in the interior not the boundary of the underlying manifold. 
%\nc

\subsection{Forward Map Discretization}\label{subsec:ForwardMapDiscrete}
Using the tools introduced in previous sections, here we propose graph-based discretizations of forward maps for elliptic and heat inverse problems on one and two-dimensional manifolds with boundary.
\paragraph{Elliptic Inverse Problem}
For the inversion problem involving the elliptic PDE in \eqref{eq:ellipticPDEdirichlet}, our goal is to learn the diffusion coefficient $\kappa$ on $\mathcal{M}$. Correspondingly, analogous to the boundary-free setting, a forward map was given by $\mathcal{F}: \theta \mapsto u$, where $\theta = \log \kappa \in (-\infty, \infty)$ and $u$ solves \eqref{eq:ellipticPDEdirichlet}. In this case, we can use the GPDM algorithm introduced in Section \ref{sec31} to obtain the discretization of the forward map. Precisely, the discretized forward map is given by 
$$
\mathcal{F}_{\epsilon, N}: \theta_N \mapsto u_N = (u_{N-B}, u_B) = \left( \tilde{\mathcal{L}}_{\epsilon,(N-B)\times(N-B)}^{-1} (f_{N-B} - \tilde{\mathcal{L}}_{\epsilon,(N-B)\times B} h_{B}), h_B \right),
$$ 
where $\theta_N = \bigl(\log(\kappa(x_1)), \ldots, \log(\kappa(x_N))\bigr)$. From this discretization process, we arrive at the relationship between the discretized posterior distribution $\mu_N^y$ and the discretized prior distribution $\mu_N$, which is given by 
$$
\frac{d\mu^y_N}{d\mu_N}(\theta_N) \propto \exp\left(-\frac{1}{2}|y-\mathcal{G}_{\epsilon, N} (\theta_N)|_{\Gamma}^2\right),
$$
where  $\mathcal{G}_{\epsilon, N}(\theta_N) = \bigl( u_N(\tilde x_1), \ldots, u_N(\tilde x_M) \bigr).$ 
For the sampling and interpolation steps, we follow the general methodology in Section \ref{subsec: setup}.

\paragraph{Heat Inversion}
 Consider first the one-dimensional heat equation with Dirichlet boundary conditions given by \eqref{eq:heatdirichlet} with $\partial \mathcal{M} = \{x_1, x_N\}$. Analogous to the prior construction, we introduce two time-independent functions $\psi_1$ and $\psi_2$ defined on $\mathcal{M}$ satisfying
\begin{align*}
\begin{cases}
\Delta_{\mathcal{M}} \psi_1 = 0,  \\
\psi_1(x_1) = 1, \\
\psi_1(x_N) = 0,
\end{cases}
\quad
\text{and}
\quad \quad
\begin{cases}
\Delta_{\mathcal{M}} \psi_2 = 0, \\
\psi_2(x_1) = 0, \\
\psi_2(x_N) = 1. 
\end{cases}
\end{align*}
Then we can obtain the solution of \eqref{eq:heatdirichlet} by solving for $w$ that satisfies,
\begin{align}\label{subsec:BIPMBheathom}
\begin{cases}
w_t(x, t) = -\Delta_\mathcal{M} w(x, t) , \quad x \in \mathcal{M}^o, t > 0,   \\
w(x_1, t) =0, ~w(x_N, t) = 0, \quad t \ge 0, \\
w(x, 0) = \theta(x)-h(x_1)\psi_1(x) - h(x_N)\psi_2(x), \quad x \in \mathcal{M}.
\end{cases}
\end{align}
To see this, notice that the function 
\begin{equation}\label{subsec:BIPMBheat1Dforwad}
u^*(x,t) \coloneqq w(x,t) + h(x_1)\psi_1(x) + h(x_N)\psi_2(x)
\end{equation}
is the solution of \eqref{eq:heatdirichlet}. 

For the two-dimensional heat equation, let $\tilde\psi_{1,\ell}$ and $\tilde\psi_{2,\ell}$ be ordered eigenfunctions of the Laplace Beltrami operator defined on $\mathcal{B}_1$ and $\mathcal{B}_2$, respectively. We can then write
\begin{align*}
h_1(x) = \sum_{\ell=1}^\infty a_\ell \tilde\psi_{1,\ell}(x), \quad x \in \mathcal{B}_1, \quad  ~\text{and}~ \quad
h_2 (x) =  \sum_{\ell=1}^\infty  b_\ell\tilde\psi_{2,\ell}(x), \quad x \in   \mathcal{B}_2,
\end{align*}
for some real coefficients $\{a_\ell\}_{\ell=1}^\infty$ and $\{b_\ell\}_{\ell=1}^\infty$. The solution to \eqref{eq:heatdirichlet} is then given by 
\begin{equation}\label{subsec:BIPMBheat2Dforwad}
u^*(x,t) \coloneqq w(x,t) +  \sum_{\ell=1}^\infty a_\ell\psi_{1,\ell}(x) + \sum_{\ell=1}^\infty b_\ell\psi_{2,\ell}(x),
\end{equation}
where $w$ is the solution of
\begin{align}\label{subsec:BIPMBhom2Dheat}
\begin{cases}
 w_t(x, t) = -\Delta_\mathcal{M} w(x, t) , \quad x \in \mathcal{M}^o, \quad t > 0,   \\
w(x, t) = 0, \quad x \in \mathcal{B}_1, \quad t \ge 0, \\w(x, t) = 0, \quad x \in \mathcal{B}_2, \quad t \ge 0, \\
w(x, 0) = \theta(x)-\sum_{\ell=1}^\infty  a_\ell\psi_{1,\ell}(x) - \sum_{\ell=1}^\infty b_\ell\psi_{2,\ell}(x), \quad x \in \mathcal{M}.
\end{cases}
\end{align}
 Here, $\psi_{1,\ell}$ is the harmonic function with boundary condition  $\tilde \psi_{1,\ell}$ in $\mathcal{B}_1$ and $0$ in $\mathcal{B}_2,$ while $\psi_{2,\ell}$ is the harmonic function with boundary condition $0$ in $\mathcal{B}_1$ and $\tilde \psi_{2,\ell}$ in $\mathcal{B}_2.$ We have hence expressed the solution of the heat equation with non-homogeneous boundary conditions as a linear superposition of basis-like functions and the solution of homogeneous heat equation. Such a decomposition will play a key role in the discretization procedures.

Now suppose that in the one-dimensional case we are given an $N$-dimensional discrete representation $\theta_N$ of the initial heat distribution of the form
$$
 \theta_N = \psi_3^{(N)} + \mu_1 \psi_1^{(N)} + \mu_2\psi_2^{(N)}, ~\text{where}~ \psi_3^{(N)} =  \frac{\sqrt{N}}{\sqrt{\sum_{n=1}^N (\tau + \lambda_n^{(N)})^{-s}}}\sum_{n=1}^N(\tau + \lambda_n^{(N)})^{-\frac{s}{2}}\zeta_n \phi_n^{(N)},
$$
which reflects our prior construction given in the previous section. To discretize the forward map of the one-dimensional heat equation, we replace the Laplace-Beltrami operator $\Delta_{\mathcal{M}}$ by $\tilde \Delta_N$ and solve for $w^{(N)} = (w_1, \ldots, w_N)$ satisfying
\begin{align}\label{interior}
\begin{cases}
\frac{\partial}{\partial t}w^{(N)} = -\tilde \Delta_N w^{(N)} ,  \\
w^{(N)}(0) = \psi_3^{(N)} , %~x \in \mathcal{M},
\end{cases}
\end{align}
This can be viewed as a discrete analog of \eqref{subsec:BIPMBheathom} as the eigenvectors of $\tilde \Delta_N$ tend to have boundary values close to zero. Given an initial condition $\psi_3^{(N)} =\theta_N - \mu_1 \psi_1^{(N)} - \mu_2\psi_2^{(N)}$, the solution of the above initial value problem can be expressed as a linear combination of the eigenvectors of $\tilde \Delta_N$. In other words, the solution at time $t$ is given by 
\begin{equation*}\label{subsubsec:BIPMBheatTrueFor}
    w^{(N)}(t) = \sum_{n = 1}^N \langle \phi_n^{(N)} ,  \psi_3^{(N)} \rangle e^{-\lambda_n^{(N)}t} \phi_n^{(N)} .
\end{equation*}
From this, and the previous observations we made in \eqref{subsec:BIPMBheat1Dforwad}, we naturally obtain the discrete approximation for the solution of \eqref{eq:heatdirichlet} given by
$$
u_N =  w^{(N)} + \mu_1\psi_1^{(N)} + \mu_2\psi_2^{(N)},
$$
which defines the discretized forward map 
$ \mathcal{F}_N: \theta_N \mapsto u_N.$
  
Using the same argument, we can discretize the forward map for the two-dimensional heat equation. Suppose we are given an $N$-dimensional discrete representation of the initial function for two-dimensional heat inverse problem, denoted by $\theta_N$, which is given by
\begin{equation}\label{2Dheatsample}
   \theta_N = \psi_3^{(N)} + \sum_{\ell=1}^L\mu_{1,\ell} \psi_{1,\ell}^{(N)} + \sum_{\ell=1}^L\mu_{2,\ell} \psi_{2,\ell}^{(N)}, ~\text{where}~ \psi_3^{(N)} =  \frac{\sqrt{N}}{\sqrt{\sum_{n=1}^N (\tau + \lambda_n^{(N)})^{-s}}}\sum_{n=1}^N(\tau + \lambda_n^{(N)})^{-\frac{s}{2}}\zeta_n \phi_n^{(N)}.  
\end{equation}
Similarly as in the one-dimensional case, with the observation \eqref{subsec:BIPMBheat2Dforwad}, the discrete approximation for the solution of \eqref{eq:heatdirichlet} in the two-dimensional case is given by 
$$
u_N  =  w^{(N)} + \sum_{\ell=1}^L\mu_{1,\ell} \psi_{1,\ell}^{(N)} + \sum_{\ell=1}^L\mu_{2,\ell} \psi_{2,\ell}^{(N)},
$$
where $w^{(N)}$ denotes the solution of \eqref{interior}.
Accordingly, the discretized forward map is given by
$
\mathcal{F}_N: \theta_N \mapsto u_N.
$
In both the one and two-dimensional cases the discretized posterior distribution $\mu_N^y$ and the discretized prior distribution $\mu_N$ satisfy the relationship
$$
\frac{d\mu^y_N}{d\mu_N}(\theta_N) \propto \exp\left(-\frac{1}{2}|y-\mathcal{G}_N(\theta_N)|_{\Gamma}^2\right),
$$
where  $\G_N(\theta_N) = \bigl( u_N(\tilde x_1), \ldots, u_N(\tilde x_M) \bigr).$ 
For the sampling and interpolation steps, we again follow the general methodology  in Section \ref{subsec: setup}.

\section{Numerical Results}\label{sec:numerics}
In this section, we provide simulation results for the numerical solution of PDE-constrained Bayesian inverse problems on a one-dimensional semi-ellipse and a two-dimensional semi-torus. For the semi-ellipse the embedding was given by  
\begin{equation}\label{1Dsemiellipse}
\iota( \alpha ) = \begin{pmatrix}
\cos\alpha \\ 3\sin\alpha
\end{pmatrix},
~\alpha \in [0, \pi],
\end{equation}
with Riemmanian metric
\begin{equation}\label{1DsemiellipseRM}
g = \sin^2(\alpha) + 9\cos^2(\alpha).
\end{equation}
The embedding for the semi-torus was given by
\begin{equation}\label{2Dsemiellipse}
\iota(\alpha, \beta) = \begin{pmatrix}
(2 + \cos\alpha)\cos\beta \\ (2 + \cos\alpha)\sin\beta  \\
\sin\alpha
\end{pmatrix},
~\alpha \in [0, 2\pi], ~\beta \in [0, \pi],
\end{equation}
where $(\alpha, \beta)$ are the intrinsic coordinates and the corresponding Riemmanian metric is given by
\begin{equation}\label{2DsemiellipseRM}
g = \begin{pmatrix}
1 & 0 \\
0 & (2 + \cos\alpha)^2
\end{pmatrix}.
\end{equation}
We use uniform grids to define point clouds over the given manifolds. For the one-dimensional semi-ellipse we used $N = 630$  points and for the two-dimensional semi-torus  $N = N_1 \times N_2 = 1296$ points in a 36 by 36 grid. The boundary of the semi-ellipse corresponds to $\alpha = 0$ or $\alpha = \pi$ and for the semi-torus the boundary corresponds to $\beta = 0$ or $\beta = \pi$. In all our numerical experiments we choose $N = M,$ that is, we assume to have observations along the entire point cloud. 

\subsection{Elliptic Inverse Problem}
For the one-dimensional elliptic inverse problem, we set the true PDE solution $u^{\dagger}$ to be
$$
u^{\dagger}(\alpha) = \sin(\alpha),
$$ 
and for the two-dimensional elliptic inverse problem, the true PDE solution was given by 
$$
 u^{\dagger}(\alpha, \beta) = 10\sin(2\alpha)\cos(\beta).
$$
According to the above true solution, the observations were given by 
$$
y_n = u^{\dagger}(\alpha_n) + \eta_n, ~\alpha_n = \iota^{-1}(x_n), ~n = 1, \ldots, N,
$$
where $\eta_n \stackrel{\text{i.i.d.}}{\sim}\mathcal{N}(0, 0.01)$. We consider several choices of true input paramter $\kappa$, and for each choice we define the right-hand side of the PDE using the  identity 
\begin{align*}
f &= \text{div}(\kappa \nabla u) =
\frac{1}{\sqrt{\text{det}g}}\partial_i\left(\kappa g^{ij}\partial_j u\sqrt{\text{det}g} \right).
\end{align*}
 Note that here and henceforth we abuse notation by referring to $u^\dagger\circ\iota$ as $u^\dagger.$
%{\color{cyan}While I understand the notation above, the convention is to write $u$ in terms of ambient coordinates $u(x)$ which means $u(x)=u\circ\iota (\iota^{-1}(x)) = u\circ\iota (\alpha)$. You can check for instance \url{http://www.math.toronto.edu/mgualt/courses/18-367/docs/DiffGeomNotes-5.pdf} that Shixiao point out to me at one instance.
%So, if $(x_1,x_2)=x=\iota(\alpha)$, then your true solution in 1D above is effectively, $u^\dagger(x) = \frac{1}{3}x_2$ which is also means, $u^\dagger\circ\iota (\alpha) = \sin(\alpha)$. Maybe we can mention that we abuse the notation here by referring to $u^\dagger\circ\iota$ as $u^\dagger$. This also applies to the entire notations in this section, as $\kappa$ is also defined in term of intrinsic coordinates.
%}

\subsubsection{One-dimensional Manifold}
For the one-dimensional elliptic inverse problem, we first augmented the given manifold point cloud data, adding $10$ ghost points at each boundary point. To model the interior term from the proposed prior, we constructed a self-tuned graph Laplacian using both the manifold point cloud and ghost points. After obtaining the graph Laplacian, we truncated it to obtain a submatrix whose indices correspond to the manifold point cloud elements. Recall that this matrix was referred to as a \emph{truncated graph Laplacian} in the previous section \ref{subsec:BIEPHeatDirichlet}. We used two nearest neighbors to construct the self-tuned graph Laplacian. For the two boundary terms, we used GPDM to obtain two harmonic functions whose boundary values were either one or zero. When constructing the weighted Laplacian matrix, we used 51 closest points in \eqref{subsec:epscalibration} and chose the value of $\epsilon$ which attained the maximum slope of $\log \bigl(T(\epsilon)\bigr) \approx \frac{1}{2}$.

We set the smoothness prior parameter to be $s = 4$, and the inverse length-scale parameter to be $\tau = 0.2$, where the choice of these values is based on empirical experiments. Hierarchical Bayesian formulations to learn these parameters could be considered \cite{sanz2020spde}. In practice, instead of using all eigenvectors of the \emph{truncated graph Laplacian}, one can use a  sufficiently large subset of eigenvectors. In our simulation study, we used 20 eigenvectors of the \emph{truncated graph Laplacian} to represent the interior term. Furthermore, to attain an acceptance rate between 40 and 60 percent, we used $\zeta = 0.01$ for the graph pCN algorithm. We computed a total of 10000 MCMC iterations with an initial 5000 burn-in period. The results are shown in Figure \ref{1}, where three different choices of diffusion coefficients are considered. 

Observe from Figure \ref{1} that the three true diffusion coefficients considered lie, for the most part, inside of the $95$ percent credible intervals. Moreover, the PDE solutions obtained using the recovered coefficients were all very close to the PDE solution with the true coefficients. To showcase the flexibility of our prior compared to the one proposed in \cite{harlim2020kernel}, additional numerical experiments were conducted. While using the same forward map approximation given by the GPDM algorithm, we employed priors proposed in \cite{harlim2020kernel} where the graph Laplacian is constructed solely from the point cloud on the manifold. We first present the recovery results for $\kappa_1(\alpha) = 2 + \cos(3\alpha)$ and $\kappa_2(\alpha) = 1 + \cos^2(\alpha)$ with the same semi-ellipse manifold as in Figure \ref{1}. All the parameter values for priors and pCN algorithms were set to be identical as before. The results are shown in Figure \ref{7}.

%We considered three different choices of diffusion coefficients. The first row of the figures below represents the recovery of the true diffusion coefficient, while the second row depicts the solution of the elliptic PDE corresponding to the recovery of the diffusion coefficient obtained in the first row. Starting from the left columns in the figures below, diffusion coefficients of the consideration were
%\begin{itemize}
%    \item $\kappa = 2 + \cos(3\alpha)$ with $\tau = 0.4$
%    \item $\kappa = 1 + \cos^2(\alpha)$ with $\tau = 0.2$
%    \item $\kappa = 1 + \alpha(\frac{\alpha-\pi}{5})$ with $\tau = 0.2$
%\end{itemize}

%\begin{figure}[H]
%\centering
%\includegraphics[width=.328\textwidth]{Elliptic1DKappa1_final.jpg}
%%\subcaption{Recovery of $\kappa = 2 + \cos(3\alpha)$}
%\includegraphics[width=.328\textwidth]{Elliptic1DKappa2_FINAL.jpg}
%%\subcaption{Recovery of $\kappa = 1 + \cos^2(\alpha)$}
%\includegraphics[width=.328\textwidth]{Elliptic1DKappa3_final.jpg}
%%\subcaption{Recovery of $\kappa = 1 + \alpha(\frac{\alpha-\pi}{5})$}
%\includegraphics[width=.328\textwidth]{Elliptic1DSol1_final.jpg}
%\includegraphics[width=.328\textwidth]{Elliptic1DSol2_final.jpg}
%\includegraphics[width=.328\textwidth]{Elliptic1DSol3_final.jpg}
%\label{1}
%\caption{One-dimensional elliptic inverse problem: Top row - Reconstruction of $\kappa$, Bottom row - Solution of an Elliptic PDE corresponding to the $\kappa$ given right above.}
%\end{figure}

\begin{figure}[htbp]
	\centering
	\begin{subfigure}{.32\textwidth}
		\centering
		\caption{Recovery of $\kappa_1(\alpha) = 2 + \cos(3\alpha)$}
		\includegraphics[width=1\textwidth]{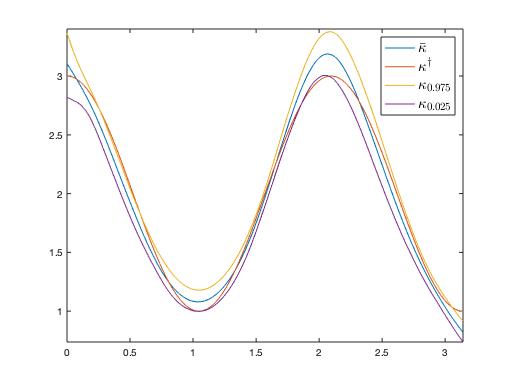}
		%\includegraphics[width=6.25cm,height=4cm]{Elliptic1DKappa1_final.jpg}
		 	%	\label{fig:i1}
	\end{subfigure}%
	%	\hfill
	\begin{subfigure}{.32\textwidth}
		\centering
		\caption{Recovery of $\kappa_2(\alpha) = 1 + \cos^2(\alpha)$}
		\includegraphics[width=1\textwidth]{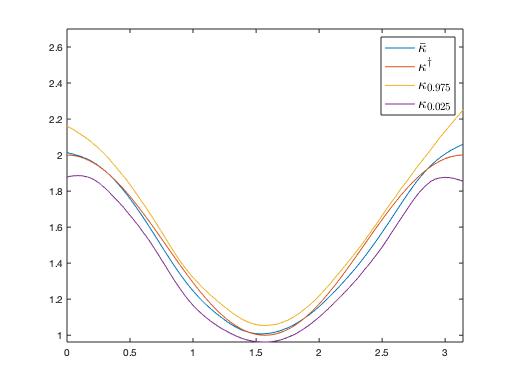}
		% 			\label{fig:i1}
	\end{subfigure}%
	%	\hfill
	\begin{subfigure}{.32\textwidth}
		\centering
		\caption{Recovery of $\kappa_3 = 1 + \frac{\alpha(\alpha - \pi)}{5}$}
		\includegraphics[width=1\textwidth]{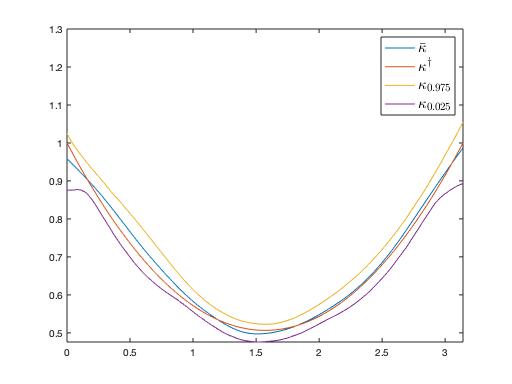}		
		% 			\label{fig:i1}
	\end{subfigure}%
	\vskip\baselineskip
	\begin{subfigure}{.32\textwidth}
		\centering
		\caption{Recovered $u_N$ of $\kappa_1$}
		\includegraphics[width=1\textwidth]{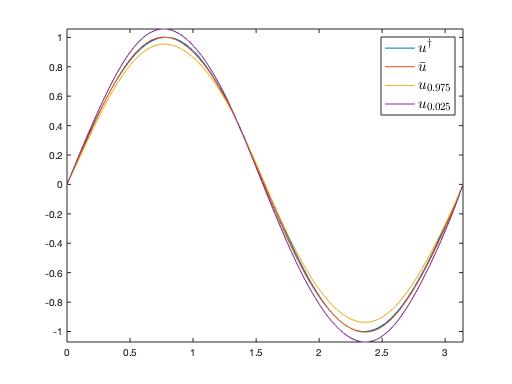}
		% 			\label{fig:i1}
	\end{subfigure}%
	%	\hfill
	\begin{subfigure}{.32\textwidth}
		\centering
				\caption{Recovered $u_N$ of $\kappa_2$}
		\includegraphics[width=1\textwidth]{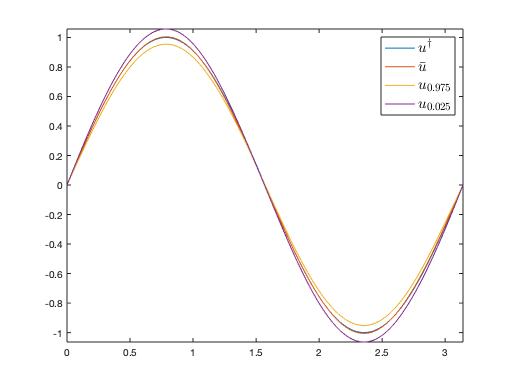}
		% 			\label{fig:i1}
	\end{subfigure}%
	\begin{subfigure}{.32\textwidth}
		\centering
		\caption{Recovered $u_N$ of $\kappa_3$}
		\includegraphics[width=1\textwidth]{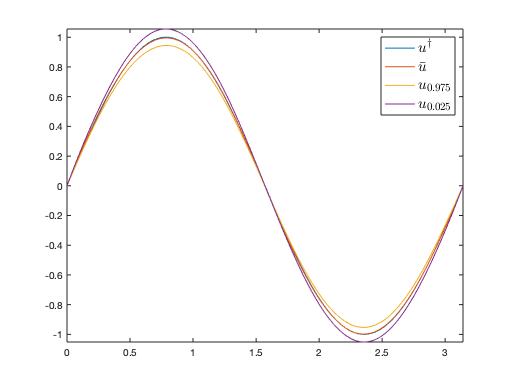}
		% 			\label{fig:i1}
	\end{subfigure}%
	\caption{One-dimensional elliptic inverse problem: Top row - Reconstruction of $\kappa$, Bottom row - Solution of an Elliptic PDE corresponding to the $\kappa$ given right above.}
	 		\label{1}
\end{figure}

%\begin{figure}[H]
%\centering
%\includegraphics[width=.328\textwidth]{Elliptic1DRuiyiKappa1_final.jpg}
%\includegraphics[width=.328\textwidth]{Elliptic1DRuiyiKappa2_final.jpg}
%\caption{One-dimensional elliptic recovery by Matern-type prior in \cite{harlim2020kernel}}
%\label{7}
%\end{figure}

\begin{figure}
	\centering
	\begin{subfigure}{.48\textwidth}
		\centering
				\caption{Recovery of $\kappa_1(\alpha) = 2 + \cos(3\alpha)$}
		\includegraphics[width=1\textwidth]{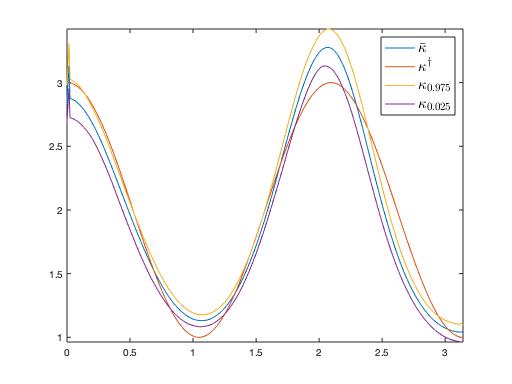}
					\label{fig:ia}
	\end{subfigure}%
	%	\hfill
	\begin{subfigure}{.48\textwidth}
		\centering
		\caption{Recovery of $\kappa_2(\alpha) = 1 + \cos^2(\alpha)$}
		\includegraphics[width=1\textwidth]{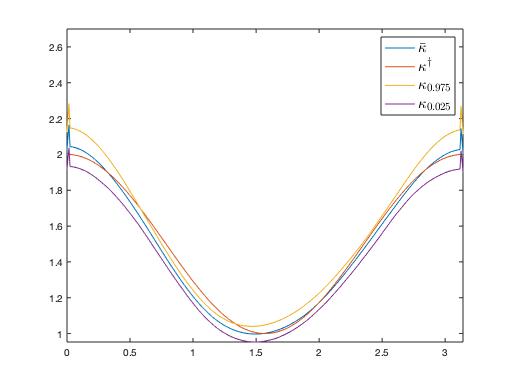}
					\label{fig:ib}
	\end{subfigure}%
	\caption{One-dimensional elliptic recovery by Mat\'ern-type prior in \cite{harlim2020kernel}.}
	 		\label{7}
\end{figure}

Figures \ref{fig:ia} and \ref{fig:ib} show that the recovered diffusion coefficients using the priors in \cite{harlim2020kernel} have artificial spikes near the boundary. The prior in \cite{harlim2020kernel} is only effective when the underlying manifold is closed, in which case the graph Laplacian used in \cite{harlim2020kernel}  approximates the Laplace Beltrami operator on the underlying manifold. However, when the underlying manifold has a boundary, the prior in \cite{harlim2020kernel}  leads to the formation of spikes near the boundaries and it does not allow for flexible modeling of boundary conditions. To illustrate this point further, we performed the Bayesian inversion procedure for the one-dimensional elliptic equation on an ellipse restricted to the first quadrant (i.e. $\alpha \in \left[0, \frac{\pi}{2}\right])$ instead of the semi-ellipse. The true diffusion coefficient was set to be $\kappa = 2 + \cos(3\alpha)$. All the parameter values were again the same as before. Figure \ref{6} demonstrates the flexibility of our proposed prior, while the reconstruction using the prior in  \cite{harlim2020kernel} exhibits artifacts near the boundary and appears to incorrectly suggest a homogeneous Neumann boundary condition for $\kappa.$

%\begin{figure}[H]
%\centering
%\includegraphics[width=.328\textwidth]{Elliptic1DKappaRuiyiCompar_final.jpg}
%\includegraphics[width=.328\textwidth]{Elliptic1DKappaCompar2_final.jpg}
%\caption{Comparison of the performance of two priors for $\kappa = 2 + \cos(3\alpha)$}
%\label{6}
%\end{figure}

\begin{figure}[ht]
	\centering
	\begin{subfigure}{.48\textwidth}
		\centering
		\includegraphics[width=1\textwidth]{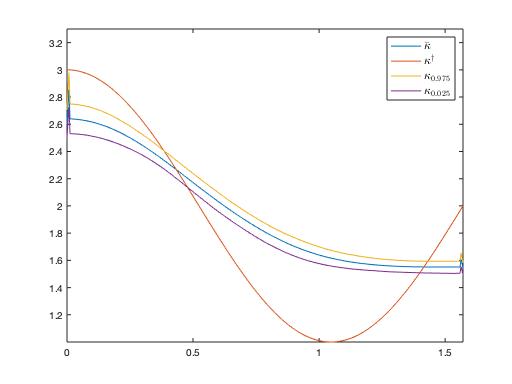}
				%\caption{Recovery of $\kappa = 2 + \cos(3\alpha).$}
		% 			\label{fig:i1}
	\end{subfigure}%
	%	\hfill
	\begin{subfigure}{.48\textwidth}
		\centering
		\includegraphics[width=1\textwidth]{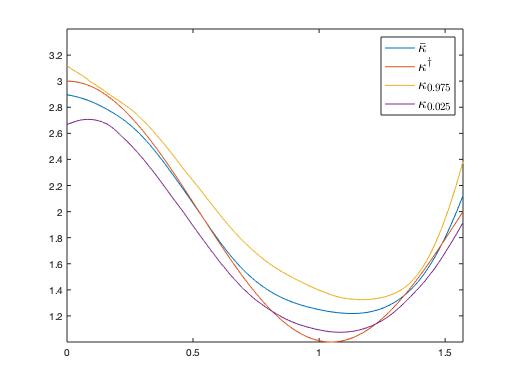}
				%\caption{Recovery of $\kappa = 1 + \cos^2(\alpha).$}
		% 			\label{fig:i1}
	\end{subfigure}%
	\caption{Comparison of the performance of two priors in the reconstruction of $\kappa_1(\alpha) = 2 + \cos(3\alpha)$. Left: using the prior proposed in \cite{harlim2020kernel}. Right: using our proposed prior.  }
	 		\label{6}
\end{figure}

%{\color{cyan} The previous width of the figures goes out of the margins; also, for aesthetic reason, I am not use to see caption below a figure so I just moved them to the top of each sub-figure.}

\subsubsection{Two-dimensional Manifold}
    Similarly, as in the one-dimensional elliptic PDE, we augmented the given point cloud with $216$ ghost points for each boundary. To define the boundary part of our prior, for each boundary (which is a circle in this example) we construct a self-tuned graph Laplacian using only the observations on the boundary with two-nearest neighbors. Then we obtain the $10$ eigenvectors corresponding to the smallest $10$ eigenvalues for each boundary. These eigenvectors will serve the role of discretized basis functions for each boundary. Solving twenty different Laplace-type of equations in \eqref{subsubsec:BIPMB2Dnonhom} using the GPDM, one can obtain $10$ basis-like functions for each boundary, which would correspond to $\{\psi_{1,\ell}^{(N)}\}_{\ell=1}^{10}$ and $\{\psi_{2,\ell}^{(N)}\}_{\ell=1}^{10}$ in  \eqref{subsubsec:2DPrior}.  When constructing the weighted Laplacian matrix, we used 128 closest points in \eqref{subsec:epscalibration} and chose the value of $\epsilon$ which attained the maximum slope of $\log \bigl(T(\epsilon) \bigr) \approx 1$.  
    
    For the interior part of the proposed prior, analogously to the one-dimensional case, we truncated a self-tuned graph Laplacian constructed using the augmented dataset and obtained a submatrix whose indices correspond to the elements in the point cloud. Taking eigenvalues and eigenvectors of this \emph{truncated graph Laplacian} would give $(\lambda_n^{(N)}, \phi_n^{(N)})$ in \eqref{subsubsec:2DPrior}. We used four-nearest neighbors to construct the self-tuned graph Laplacian. We set the smoothness prior parameter to be $s = 4$, and the inverse length-scale parameter to be $\tau = 0.24$, where the choice of these values is based on empirical experiments. To attain an acceptance rate between 40 and 60 percent, we used $\zeta = 0.001$ for the graph pCN algorithm. A total of $150000$ MCMC iterations with initial $75000$ burn-in iterations were run to obtain our results. 
    
    Starting from the top row of Figure \ref{5}, the panels in the first row  represent the true diffusion coefficient and the solution of the elliptic PDE based on the true diffusion coefficient. The second row represents, from left to right, the posterior mean of the MCMC samples, 2.5 percentile of MCMC samples and 97.5 percentile of MCMC samples. In other words, the second row portrays the credible interval for the true diffusion coefficient function. The third row corresponds to the solution of the elliptic PDE equation based on the diffusion coefficients given in the second row. The last two plots in the fourth row depict the difference between the true diffusion coefficient and the posterior mean of the MCMC samples, and the difference between the true solution with the one based on the posterior mean of the MCMC samples.

    %The last row, from left to right, corresponds to the posterior mean, 2.5th percentile and 97.5 percentile of the all solutions corresponding to the MCMC samples of diffusion coefficients. 
    
    For the two-dimensional elliptic inverse problem, there is a larger error in the recovery of the true diffusion coefficient. This is unsurprising, since as one can see from the solutions corresponding to the MCMC samples of the diffusion coefficients, one can find several different diffusion coefficients that lead to approximate solutions that are close to the solution of the elliptic PDE, which is a manifestation of the ill-posedness of this elliptic inverse problem.

%\begin{figure}[H]
%\centering
%\includegraphics[width=.33\textwidth]{Elliptic2DTrueKappa.jpg}
%\includegraphics[width=.33\textwidth]{Elliptic2DSol.jpg}
%\includegraphics[width=.328\textwidth]{Elliptic2DKappaPM.jpg}
%\includegraphics[width=.328\textwidth]{Elliptic2DKappaPML.jpg}
%\includegraphics[width=.328\textwidth]{Elliptic2DKappaPMU.jpg}
%\includegraphics[width=.328\textwidth]{Elliptic2DKappaPMsol.jpg}
%\includegraphics[width=.328\textwidth]{Elliptic2DKappaPMLsol.jpg}
%\includegraphics[width=.328\textwidth]{Elliptic2DKappaPMUsol.jpg}
%\includegraphics[width=.328\textwidth]{Elliptic2DUM.jpg}
%\includegraphics[width=.328\textwidth]{Elliptic2DUML.jpg}
%\includegraphics[width=.328\textwidth]{Elliptic2DUMU.jpg}
%%\includegraphics[width=.33\textwidth]{Elliptic2DlogKappa-PM.jpg}
%%\includegraphics[width=.33\textwidth]{Elliptic2DU2-UM.jpg}
%\caption{Two-dimensional elliptic inverse problem}
%\label{5}
%\end{figure}

\begin{figure}[htbp]
	\centering
	\begin{subfigure}{.32\textwidth}
		\centering
		\caption{True diffusion coefficient $\kappa$}
		\includegraphics[width=1\textwidth,height=4cm]{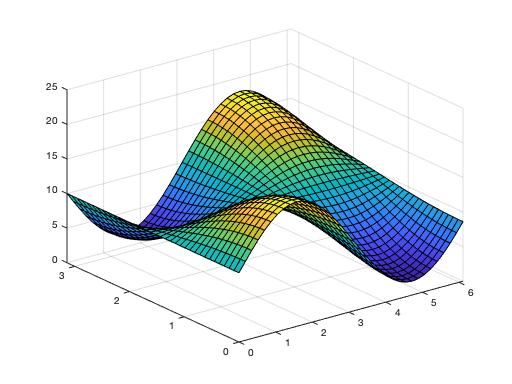}
		 	%	\label{fig:i1}
	\end{subfigure}%
	%	\hfill
	\begin{subfigure}{.32\textwidth}
		\centering
				\caption{True PDE solution $u$}
		\includegraphics[width=1\textwidth,height=4cm]{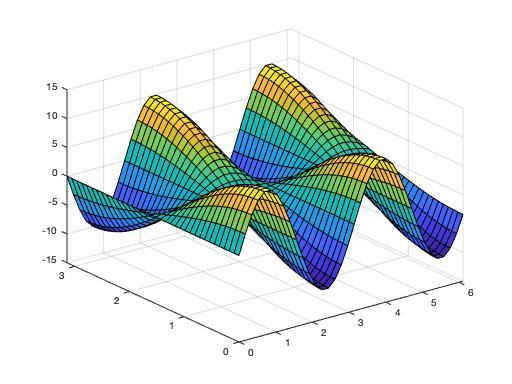}
		% 			\label{fig:i1}
	\end{subfigure}%
	\vskip\baselineskip
	%	\hfill
	\begin{subfigure}{.32\textwidth}
		\centering
				\caption{Posterior mean for $\kappa$}
		\includegraphics[width=1\textwidth,height=4cm]{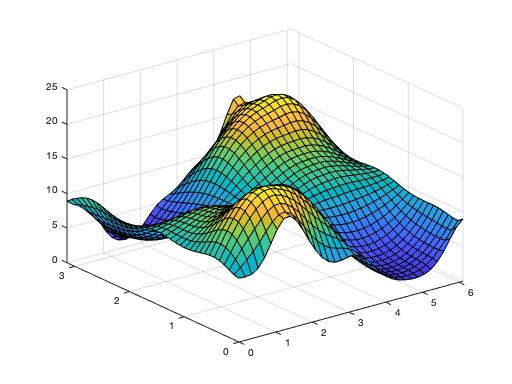}
					\label{figc}
	\end{subfigure}%
	\begin{subfigure}{.32\textwidth}
		\centering
		\caption{$2.5\%$ post.~percentile for $\kappa$}
		\includegraphics[width=1\textwidth,height=4cm]{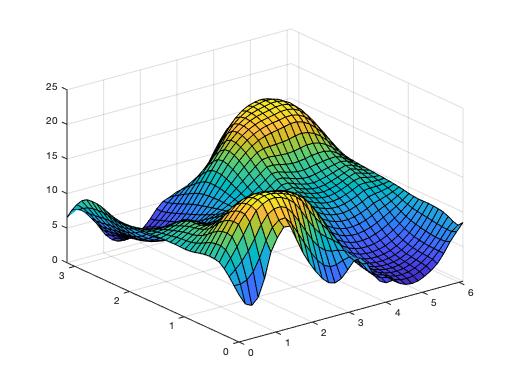}
				\label{figd}
	\end{subfigure}%
	%	\hfill
	\begin{subfigure}{.32\textwidth}
		\centering
				\caption{$97.5\%$ post.~percentile for $\kappa$}
		\includegraphics[width=1\textwidth,height=4cm]{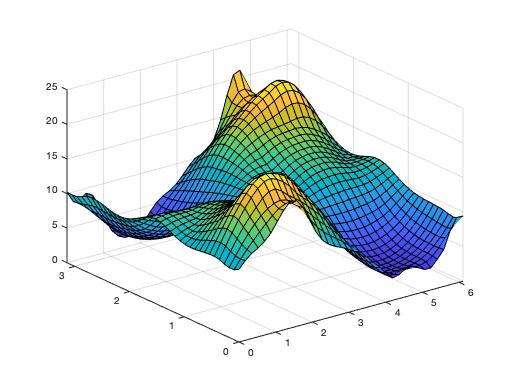}
				\label{fige}
	\end{subfigure}%
	\vskip\baselineskip
	\begin{subfigure}{.32\textwidth}
		\centering
				\caption{$u_N$ with $\kappa$ in Figure \ref{figc}}
		\includegraphics[width=1\textwidth,height=4cm]{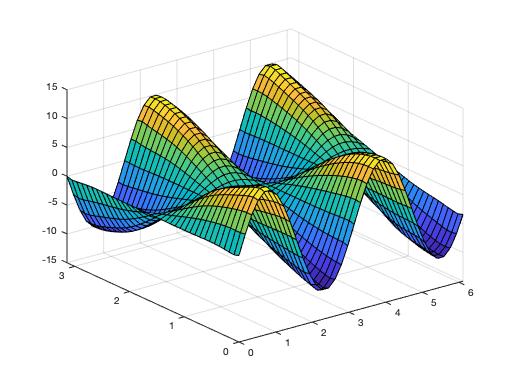}
		% 			\label{fig:i1}
	\end{subfigure}%
	\begin{subfigure}{.32\textwidth}
		\centering
				\caption{$u_N$ with $\kappa$ in Figure \ref{figd}}
		\includegraphics[width=1\textwidth,height=4cm]{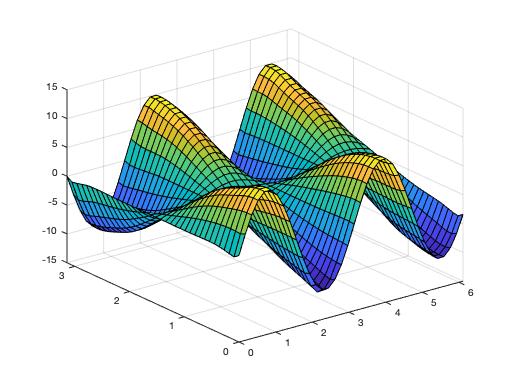}
		% 			\label{fig:i1}
	\end{subfigure}%
	%	\hfill
	\begin{subfigure}{.32\textwidth}
		\centering
				\caption{$u_N$ with $\kappa$ in Figure \ref{fige}}
		\includegraphics[width=1\textwidth,height=4cm]{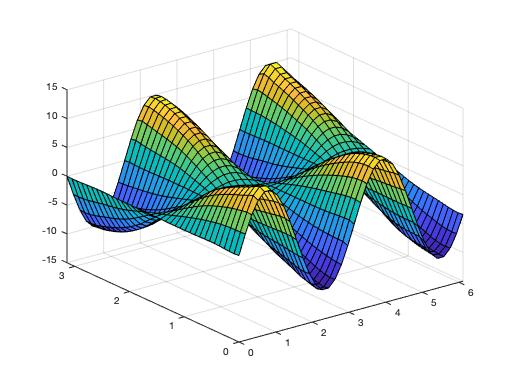}
		% 			\label{fig:i1}
	\end{subfigure}%
	\vskip\baselineskip
	\begin{subfigure}{.4\textwidth}
		\centering
				\caption{Error in diffusion coefficient}
		\includegraphics[width=1\textwidth,height=4cm]{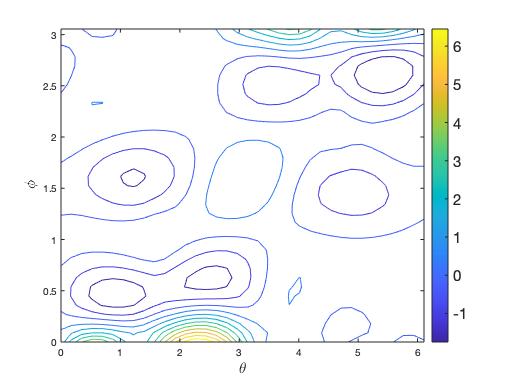}
		 			%\label{fig:7j}
	\end{subfigure}%
	\begin{subfigure}{.4\textwidth}
		\centering
				\caption{Error in PDE solution}
		\includegraphics[width=1\textwidth,height=4cm]{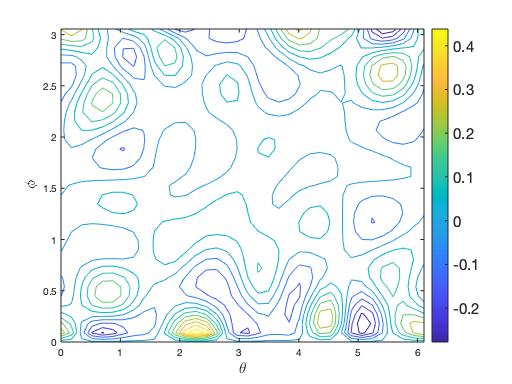}
		 			%\label{fig:7k}
	\end{subfigure}%
	%\begin{subfigure}{.32\textwidth}
	%	\centering
	%			\caption{Posterior mean for $u_N$}
	%	\includegraphics[width=1\textwidth,height=4cm]{ELLIPTIC2D_SOLUMM_FINAL.jpg}
		% 			\label{fig:i1}
	%\end{subfigure}%
	%	\begin{subfigure}{.32\textwidth}
	%	\centering
	%			\caption{$2.5\%$ post.~percentile for $u_N$}
	%	\includegraphics[width=1\textwidth,height=4cm]{ELLIPTIC2D_UMML_FINAL.jpg}
		% 			\label{fig:i1}
	%\end{subfigure}%
	%	\begin{subfigure}{.32\textwidth}
	%	\centering
    %	\caption{$97.5\%$ post.~percentile for $u_N$}
	%	\includegraphics[width=1\textwidth,height=4cm]{ELLIPTIC2D_UMMU_FINAL.jpg}
		% 			\label{fig:i1}
	%\end{subfigure}%
	\caption{Elliptic inverse problem on a semi-torus for $\kappa(\alpha, \beta) = 10 + 8\sin(\alpha)\cos(\beta)$.}
	 		\label{5}
\end{figure}

\subsection{Heat Inverse Problem}
For the numerical simulations of heat inversion, given an initial heat function $u_0$, observations were obtained by 
$$
y_n =  \tilde u  (\alpha_n, t^*) + \eta_n, ~\alpha_n  = \iota^{-1}(x_n), ~n = 1, \ldots,  N ,
$$
where $\eta_n \stackrel{\text{i.i.d.}}{\sim}\mathcal{N}(0, 0.01)$ and $\tilde u$ represents an approximate solution of the heat equation given in \eqref{eq:heatdirichlet}. For the one-dimensional semi-ellipse, we used the following explicit formula to compute the approximate solution over the point cloud at time $t^* > 0$, denoted by $\tilde u_N = (\tilde u(\alpha_1, t^*), \cdots, \tilde u(\alpha_N, t^*))$, 
$$
\tilde u_N \coloneqq \sum_{n = 1}^N \langle \phi_n^{(N)} ,  \psi_3^{(N)} \rangle e^{-\lambda_n^{(N)}t^*} \phi_n^{(N)} + u_0(\alpha_1)\psi_1^{(N)} + u_0(\alpha_N)\psi_2^{(N)},
$$
where $(\lambda_n^{(N)}, \phi_n^{(N)})_{n=1}^N$ are eigenpairs of the \emph{truncated graph Laplacian} constructed using point cloud $\{x_n\}_{n=1}^N$ with two boundary points $\{x_1, x_N\}$. The definition of $\psi_1^{(N)}$, $\psi_2^{(N)}$ and $\psi_3^{(N)}$ are the same as in the Section \ref{subsec:BIEPHeatDirichlet}. Similarly, for the two-dimensional semi-torus, the approximate solution $\tilde u_N$ over the point cloud at time $t^* > 0$, denoted by $\tilde u_N = (\tilde u(\alpha_1, \beta_1, t^*), \cdots, \tilde u(\alpha_{N_1}, \beta_{N_2}, t^*))$ with $N = N_1\times N_2$ is given by 
$$
\tilde u_N \coloneqq  \sum_{n = 1}^N \langle \phi_n^{(N)} ,  \psi_3^{(N)} \rangle e^{-\lambda_n^{(N)}t^*} \phi_n^{(N)} + \sum_{\ell=1}^L\tilde\mu_{1,\ell} \psi_{1,\ell}^{(N)} + \sum_{\ell=1}^L\tilde\mu_{2,\ell} \psi_{2,\ell}^{(N)},
$$
where the leftmost term in the right hand side is analogously defined as in the one-dimensional manifold setting given above. The coefficients $\{\tilde\mu_{1,\ell}\}_{l=1}^L$ and $\{\tilde\mu_{1,\ell}\}_{l=1}^L$ were obtained by regressing $u_0^{(N)} = (u_0(x_1), \cdots, u_0(x_N))$ on $\{\phi_n^{(N)}\}_{n=1}^N, \{\psi_{1,l}^{(N)}\}_{l=1}^L$ and $\{\psi_{2,l}^{(N)}\}_{l=1}^L$. \nc
In our numerical experiments, we constructed $10$ basis-like functions for each boundary of the two-dimensional semi-torus, i.e., $L = 10$. For our implementation, we used $20$ eigenvectors of the \emph{truncated graph Laplacian}. 

\subsubsection{One-dimensional Manifold}
The prior construction was exactly the same as in the one-dimensional elliptic problem. We set the smoothness prior parameter to be $s = 6$, and the inverse length-scale parameter to be $\tau = 0.3$, where the choice of these values is based on empirical experiments. Furthermore, to attain an MCMC acceptance rate of roughly 50 percent, we used $\zeta = 0.005$ for the graph pCN algorithm. A total of 20000 MCMC iterations with initial 10000 burn-in iterations were run to obtain our results. For the heat equation, the degree of ill-posedness of the inverse problem is closely related to the time $t^* > 0$ at which we observe the data. The smaller the time $t^* > 0 $ is, the easier the inversion.

We considered three different choices of initial heat functions. Starting from the left-most column of  Figure \ref{2}, plots in each column respectively represent initial heat function with observed data, true initial function with the posterior mean/2.5th and 97.5th percentile of MCMC samples, and the true solution with the solutions corresponding to the posterior mean, 2.5th and 97.5 percentile of MCMC samples. Specifically, the leftmost column in  Figure \ref{2} represents the initial heat function $u_0$ with its corresponding noise-free observation $u$ and noisy observation $y$. The middle column represents the true initial function $u_0^{\dagger}$ with the posterior mean estimate $\bar u_0$, 2.5th percentile $u_0^{0.025}$, and 97.5th percentile $u_0^{0.975}$ of MCMC samples. The rightmost column portrays the approximate solution based on the true initial heat function at time $t$, denoted by $u_t^{\dagger}$ with approximate solutions based on the posterior mean estimate, 2.5th percentile and 97.5th percentile, respectively denoted by $\bar u_t$, $u_t^{0.025}$ and $u_t^{0.975}$. 
%\begin{figure}[H]
%\centering
%\includegraphics[width=.328\textwidth]{Heat1DData1_final.jpg}
%\includegraphics[width=.328\textwidth]{Heat1DInit1_final.jpg}
%\includegraphics[width=.328\textwidth]{Heat1DSol1_final.jpg}
%\includegraphics[width=.328\textwidth]{Heat1DData2_final.jpg}
%\includegraphics[width=.328\textwidth]{Heat1DInit2_final.jpg}
%\includegraphics[width=.328\textwidth]{Heat1DSol2_final.jpg}
%\includegraphics[width=.328\textwidth]{Heat1DData3_final.jpg}
%\includegraphics[width=.328\textwidth]{Heat1DInit3_final.jpg}
%\includegraphics[width=.328\textwidth]{Heat1DSol3_final.jpg}
%\label{2}
%\caption{One-dimensional heat inverse problem}
%\end{figure}

\begin{figure}[htbp]
	\centering
	\begin{subfigure}{.32\textwidth}
		\centering
		\caption{$u_0 = 10\sin(\alpha)+2$}
		\includegraphics[width=1\textwidth,height=4cm]{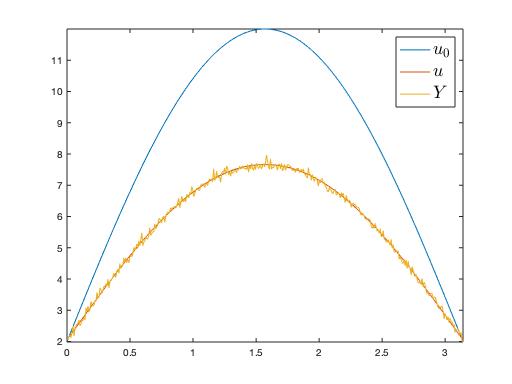}
		 	%	\label{fig:i1}
	\end{subfigure}%
	%	\hfill
	\begin{subfigure}{.32\textwidth}
		\centering
				\caption{Recovery of $u_0 = 10\sin(\alpha)+2$}
		\includegraphics[width=1\textwidth,height=4cm]{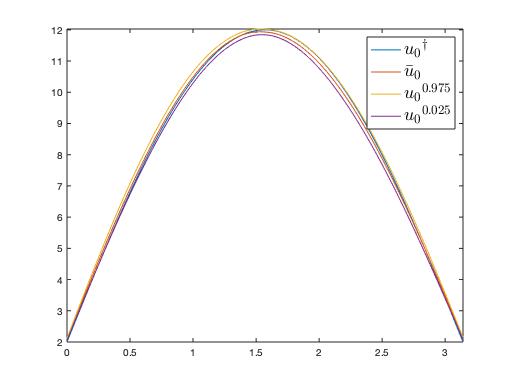}
		% 			\label{fig:i1}
	\end{subfigure}%
	%	\hfill
	\begin{subfigure}{.32\textwidth}
		\centering
				\caption{Solution $u(t)$ at $t = 10$}
		\includegraphics[width=1\textwidth,height=4cm]{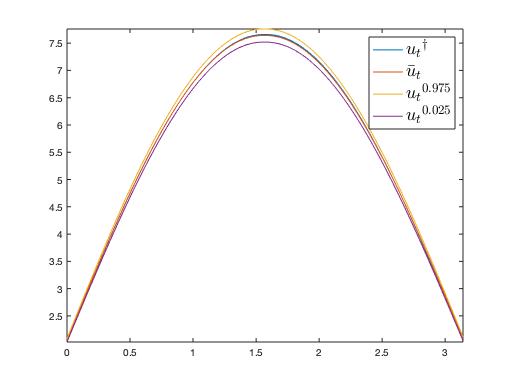}
		% 			\label{fig:i1}
	\end{subfigure}%
	\vskip\baselineskip
	\begin{subfigure}{.32\textwidth}
		\centering
				\caption{$u_0 = 10\sin(2\alpha)+2$}
		\includegraphics[width=1\textwidth,height=4cm]{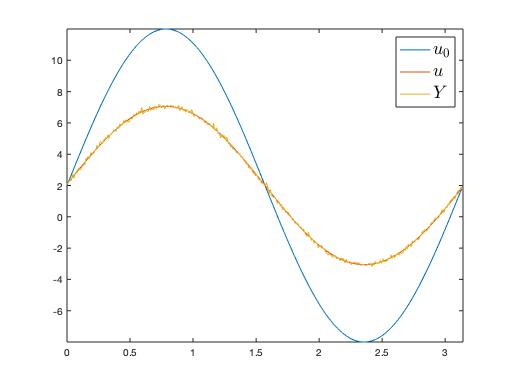}
		% 			\label{fig:i1}
	\end{subfigure}%
	%	\hfill
	\begin{subfigure}{.32\textwidth}
		\centering
				\caption{Recovery of $u_0 = 10\sin(2\alpha)+2$}
		\includegraphics[width=1\textwidth,height=4cm]{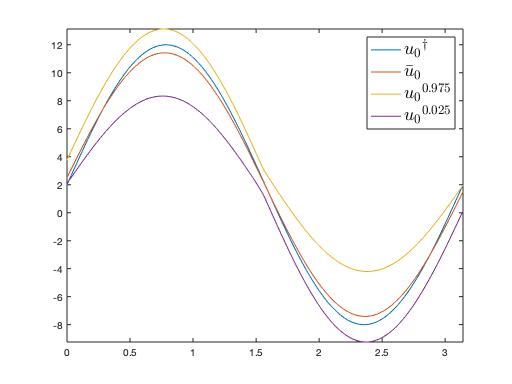}
		% 			\label{fig:i1}
	\end{subfigure}%
	\begin{subfigure}{.32\textwidth}
		\centering
				\caption{Solution $u(t)$ at $t = 3$}
		\includegraphics[width=1\textwidth,height=4cm]{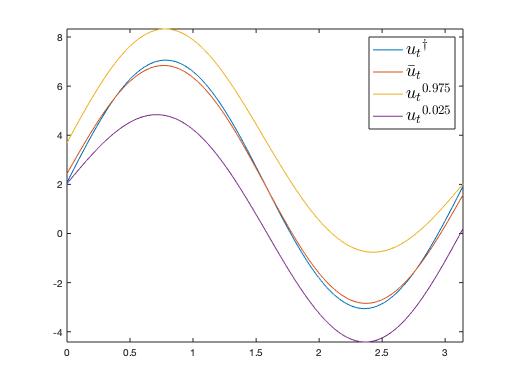}
		% 			\label{fig:i1}
	\end{subfigure}%
		\vskip\baselineskip
	\begin{subfigure}{.32\textwidth}
		\centering
				\caption{$u_0 = 10\cos(\alpha)+2$}
		\includegraphics[width=1\textwidth,height=4cm]{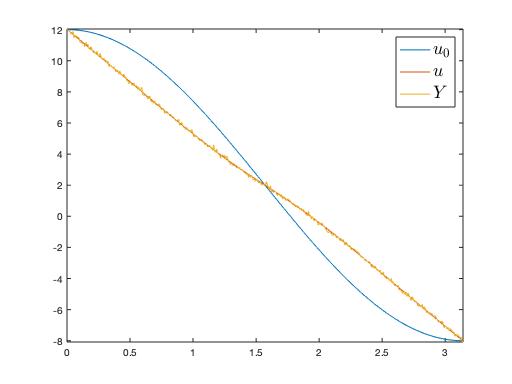}
		% 			\label{fig:i1}
	\end{subfigure}%
	%	\hfill
	\begin{subfigure}{.32\textwidth}
		\centering
				\caption{Recovery of $u_0 = 10\cos(\alpha)+2$}
		\includegraphics[width=1\textwidth,height=4cm]{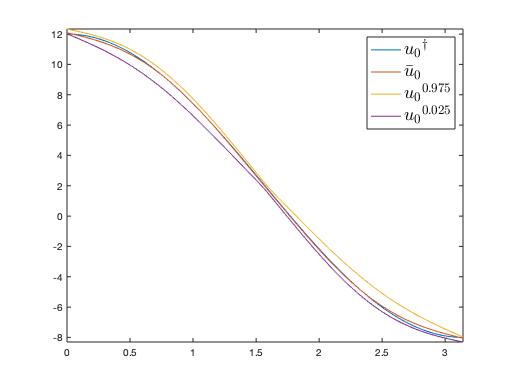}
		% 			\label{fig:i1}
	\end{subfigure}%
	\begin{subfigure}{.32\textwidth}
		\centering
				\caption{Solution $u(t)$ at $t = 5$}
		\includegraphics[width=1\textwidth,height=4cm]{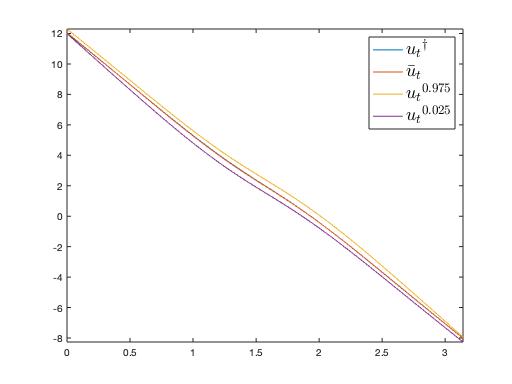}
		% 			\label{fig:i1}
	\end{subfigure}%
	\caption{Heat inverse problem on a semi-ellipse.}
	 		\label{2}
\end{figure}

Figure \ref{2} shows that all of the true initial heat functions were captured inside of the 95 percent credible intervals for all three cases. Moreover, the corresponding solutions were very close to the true solution.

\subsubsection{Two-dimensional Manifold}
The practical implementation of the two-dimensional prior was analogous to that of the two-dimensional elliptic inverse problem. We set the smoothness  prior  parameter to be $s = 4$, and the inverse length-scale parameter to be $\tau = 0.3 ~\text{or}~ 0.012$, where the choice of these values is based on empirical experiments. Furthermore, to attain an MCMC acceptance rate of roughly $50$ percent, we used $\zeta = 0.006$ for the graph pCN algorithm. A total of $100000$ number of MCMC iterations with initial $50000$ burn-in iterations were run to obtain these results. 

Starting from the top row of Figures \ref{3} and \ref{4}, the panels in the first row  represent, from left to right, the true initial heat function, approximate solution of the heat equation based on the true heat initial function, noise-incorporated approximate solution, i.e. observation data. The second row represents from left to right, posterior mean of the MCMC samples, $2.5$ percentile of MCMC samples, and $97.5$ percentile of MCMC samples. In other words, the second row portrays the credible interval for the true initial heat function. The third row corresponds to the approximate solution of the heat equation based on the initial heat functions given in the second row. Finally, the last two plots in the fourth row depict the difference between the true initial function and the posterior mean of the MCMC samples, and the difference between the approximate solution based on the true initial function with the one based on the posterior mean of the MCMC samples. Figures \ref{3} and Figure \ref{4} show that our proposed methodology led to reasonably successful recovery of the parameter of interest.

\begin{figure}[htbp]
	\centering
	\begin{subfigure}{.32\textwidth}
		\centering
				\caption{True initial heat $u_0$}
		\includegraphics[width=1\textwidth,height=4cm]{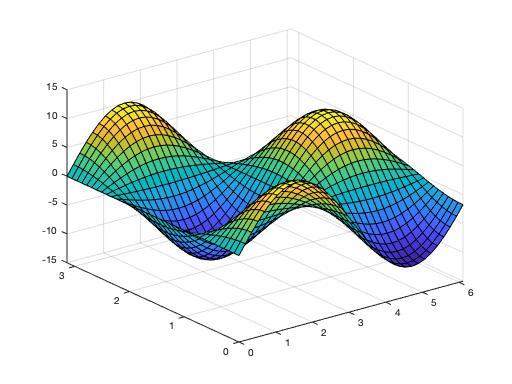}
		 		\label{fig:7a}
	\end{subfigure}%
	%	\hfill
	\begin{subfigure}{.32\textwidth}
		\centering
				\caption{PDE solution $u(t^*)$ based on $u_0$}
		\includegraphics[width=1\textwidth,height=4cm]{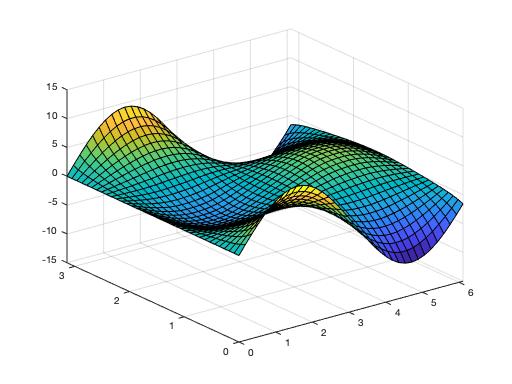}
		 			\label{fig:7b}
	\end{subfigure}%
	%	\hfill
	\begin{subfigure}{.32\textwidth}
		\centering
				\caption{Observed data $y$}
		\includegraphics[width=1\textwidth,height=4cm]{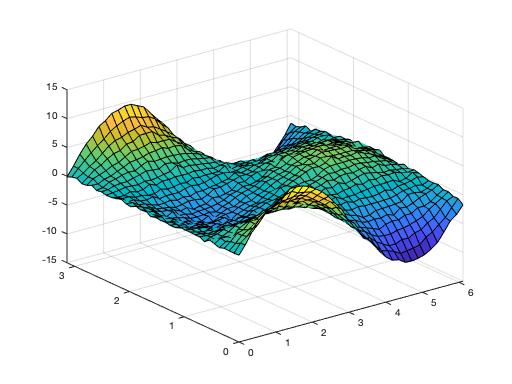}
		 			\label{fig:7c}
	\end{subfigure}%
	\vskip\baselineskip
	\begin{subfigure}{.32\textwidth}
		\centering
				\caption{Posterior mean for $u_0$}
		\includegraphics[width=1\textwidth,height=4cm]{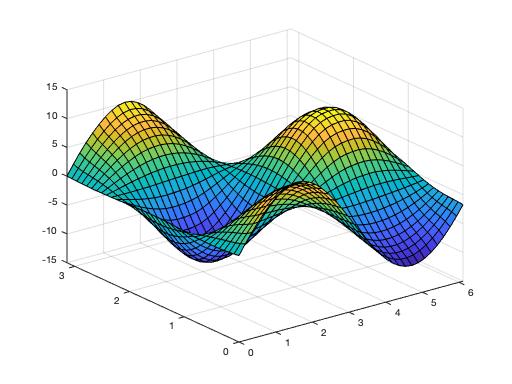}
		 			\label{fig:7d}
	\end{subfigure}%
	%	\hfill
	\begin{subfigure}{.32\textwidth}
		\centering
				\caption{$2.5\%$ post. percentile for $u_0$}
		\includegraphics[width=1\textwidth,height=4cm]{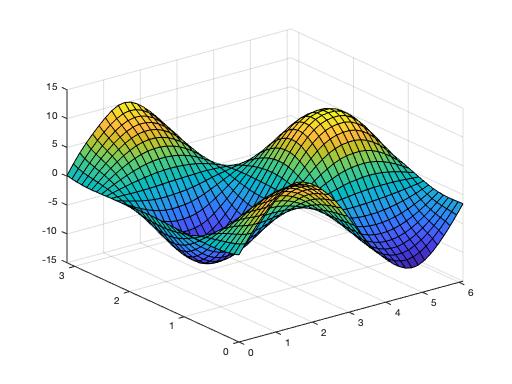}
		 			\label{fig:7e}
	\end{subfigure}%
	\begin{subfigure}{.32\textwidth}
		\centering
				\caption{$97.5\%$ post. percentile for $u_0$}
		\includegraphics[width=1\textwidth,height=4cm]{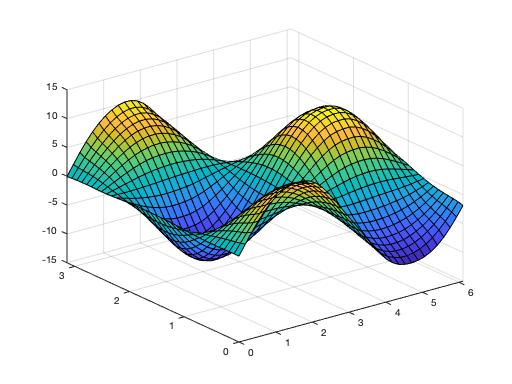}
		 			\label{fig:7f}
	\end{subfigure}%
	\vskip\baselineskip
	\begin{subfigure}{.32\textwidth}
		\centering
				\caption{$u_N(t^*)$ with $u_0$ in Figure \ref{fig:7d}}
		\includegraphics[width=1\textwidth,height=4cm]{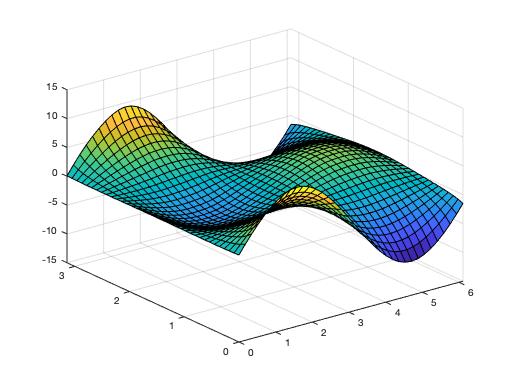}
		 			\label{fig:7g}
	\end{subfigure}%
	%	\hfill
	\begin{subfigure}{.32\textwidth}
		\centering
				\caption{$u_N(t^*)$ with $u_0$ in Figure \ref{fig:7e}}
		\includegraphics[width=1\textwidth,height=4cm]{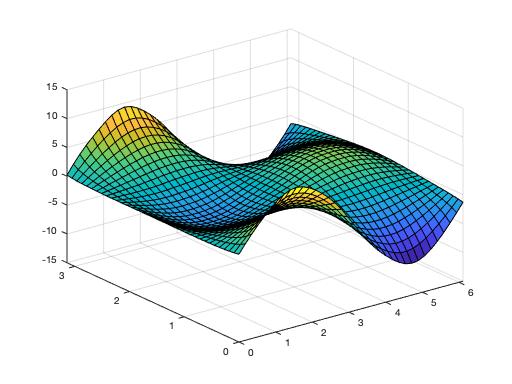}
		 			\label{fig:7h}
	\end{subfigure}%
	\begin{subfigure}{.32\textwidth}
		\centering
				\caption{$u_N(t^*)$ with $u_0$ in Figure \ref{fig:7f}}
		\includegraphics[width=1\textwidth,height=4cm]{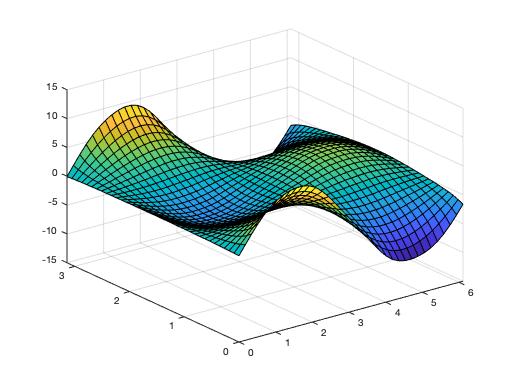}
		 			\label{fig:7i}
	\end{subfigure}%
	\vskip\baselineskip
		\begin{subfigure}{.4\textwidth}
		\centering
				\caption{Error in initial condition}
		\includegraphics[width=1\textwidth,height=4cm]{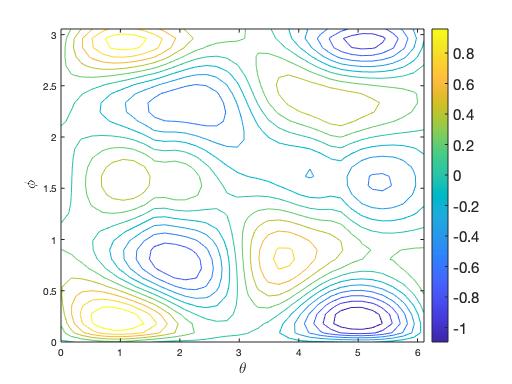}
		 			\label{fig:7j}
	\end{subfigure}%
	\begin{subfigure}{.4\textwidth}
		\centering
				\caption{Error in PDE solution}
		\includegraphics[width=1\textwidth,height=4cm]{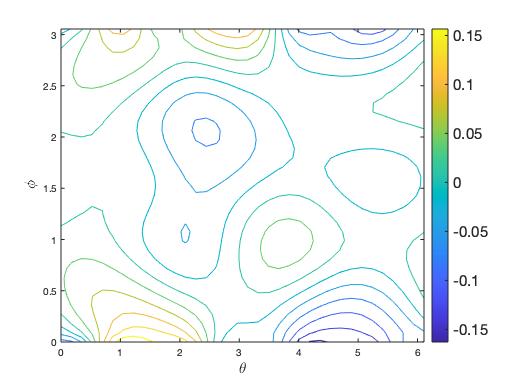}
		 			\label{fig:7k}
	\end{subfigure}%
	\vskip\baselineskip
	\caption{Heat inverse problem on a semi-torus:  $u_0(\alpha, \beta) = 10\sin(\alpha)\cos(2\beta)$ with $\tau = 0.3$, $t^* = 5$.  }
	 		\label{3}
\end{figure}

%\begin{figure}[H]
%\centering
%\includegraphics[width=.328\textwidth]{HEAT2DU0(2+sincos).jpg}
%\includegraphics[width=.328\textwidth]{HEAT2DUt(2+sincos).jpg}
%\includegraphics[width=.328\textwidth]{HEAT2DDATA(2+sincos).jpg}
%\includegraphics[width=.328\textwidth]{HEAT2DUM(2+sincos).jpg}
%\includegraphics[width=.328\textwidth]{HEAT2DU0MU(2+sincos).jpg}
%\includegraphics[width=.328\textwidth]{HEAT2DU0ML(2+sincos).jpg}
%\includegraphics[width=.328\textwidth]{HEAT2DUTM(2+sincos).jpg}
%\includegraphics[width=.328\textwidth]{HEAT2DUTMU(2+sincos).jpg}
%\includegraphics[width=.328\textwidth]{HEAT2DUTML(2+sincos).jpg}
%\includegraphics[width=.328\textwidth]{HEAT2DU0-U0M(2+sincos).jpg}
%\includegraphics[width=.328\textwidth]{HEAT2DUT-UTM(2+sincos).jpg}
%\caption{Two-dimensional Heat Inverse Problem: $u_0(\theta, \phi) = 2 + \sin(\theta)\cos(\phi)$}
%\label{4}
%\end{figure}

\begin{figure}[htbp]
	\centering
	\begin{subfigure}{.32\textwidth}
		\centering
				\caption{True initial heat $u_0$}
		\includegraphics[width=1\textwidth,height=4cm]{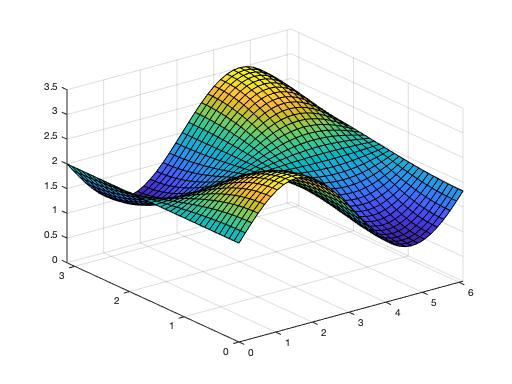}
		 		\label{fig:8a}
	\end{subfigure}%
	%	\hfill
	\begin{subfigure}{.32\textwidth}
		\centering
				\caption{PDE solution $u(t^*)$ based on $u_0$}
		\includegraphics[width=1\textwidth,height=4cm]{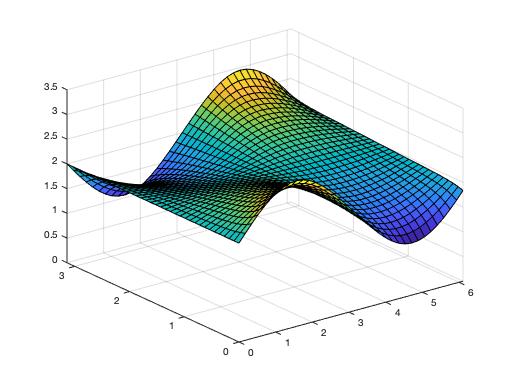}
		 			\label{fig:8b}
	\end{subfigure}%
	%	\hfill
	\begin{subfigure}{.32\textwidth}
		\centering
				\caption{Observed data $y$}
		\includegraphics[width=1\textwidth,height=4cm]{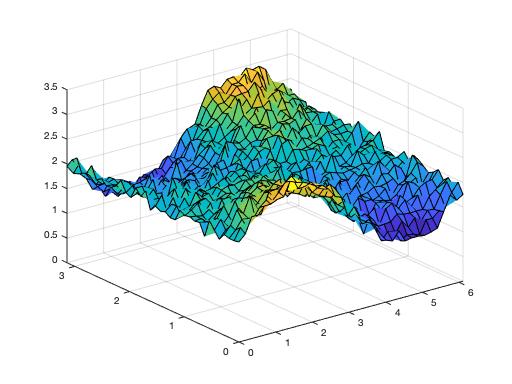}
		 			\label{fig:8c}
	\end{subfigure}%
	\vskip\baselineskip
	\begin{subfigure}{.32\textwidth}
		\centering
				\caption{Posterior mean for $u_0$}
		\includegraphics[width=1\textwidth,height=4cm]{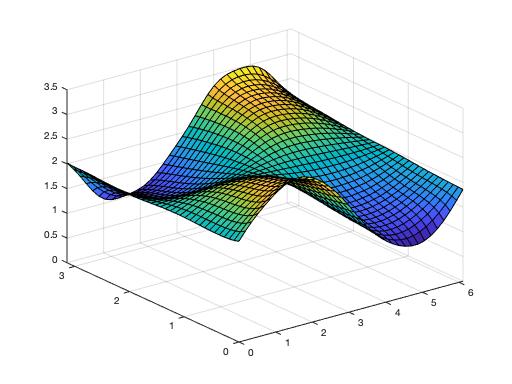}
		 			\label{fig:8d}
	\end{subfigure}%
	%	\hfill
	\begin{subfigure}{.32\textwidth}
		\centering
				\caption{$2.5\%$ post. percentile for $u_0$}
		\includegraphics[width=1\textwidth,height=4cm]{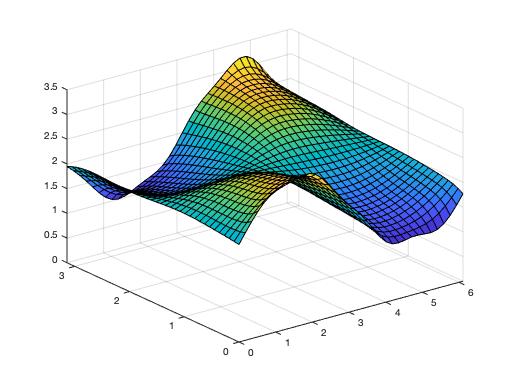}
		 			\label{fig:8e}
	\end{subfigure}%
	\begin{subfigure}{.32\textwidth}
		\centering
				\caption{$97.5\%$ post. percentile for $u_0$}
		\includegraphics[width=1\textwidth,height=4cm]{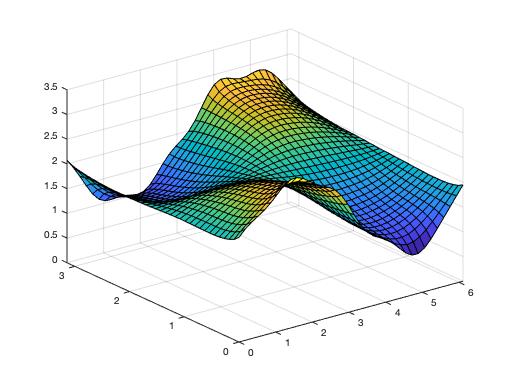}
		 			\label{fig:8f}
	\end{subfigure}%
	\vskip\baselineskip
	\begin{subfigure}{.32\textwidth}
		\centering
				\caption{$u_N(t^*)$ with $u_0$ in Figure \ref{fig:8d}}
		\includegraphics[width=1\textwidth,height=4cm]{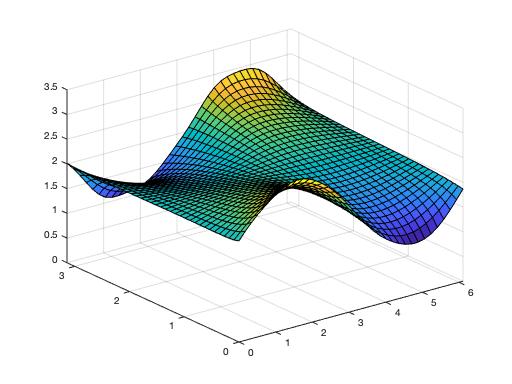}
		 			\label{fig:8g}
	\end{subfigure}%
	%	\hfill
	\begin{subfigure}{.32\textwidth}
		\centering
				\caption{$u_N(t^*)$ with $u_0$ in Figure \ref{fig:8e}}
		\includegraphics[width=1\textwidth,height=4cm]{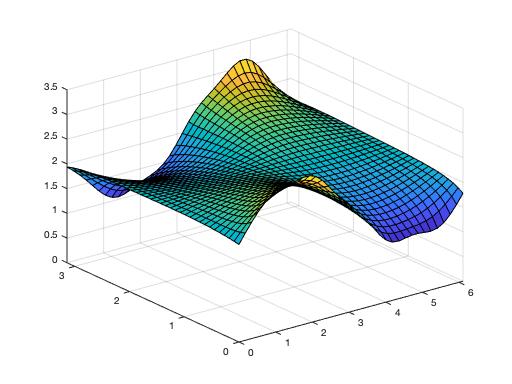}
		 			\label{fig:8h}
	\end{subfigure}%
	\begin{subfigure}{.32\textwidth}
		\centering
				\caption{$u_N(t^*)$ with $u_0$ in Figure \ref{fig:8f}}
		\includegraphics[width=1\textwidth,height=4cm]{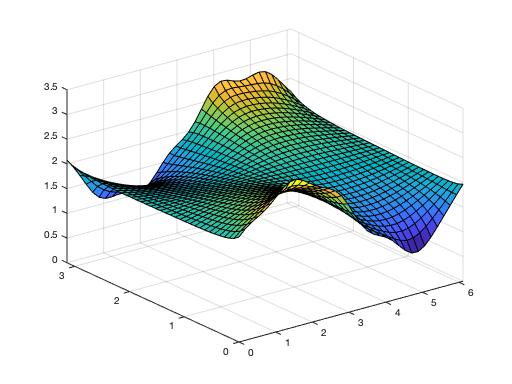}
		 			\label{fig:8i}
	\end{subfigure}%
	\vskip\baselineskip
		\begin{subfigure}{.4\textwidth}
		\centering
			\caption{Error in initial condition}
		\includegraphics[width=1\textwidth,height=4cm]{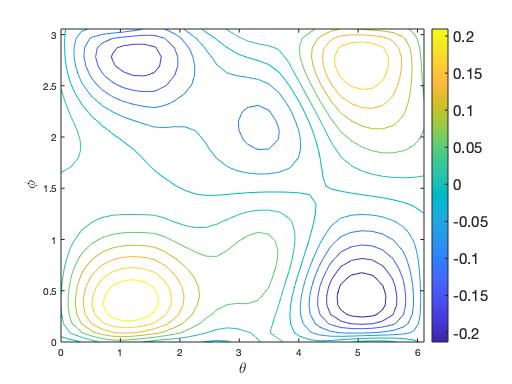}
		 			\label{fig:8j}
	\end{subfigure}%
	\begin{subfigure}{.4\textwidth}
		\centering
				\caption{Error in PDE solution}
		\includegraphics[width=1\textwidth,height=4cm]{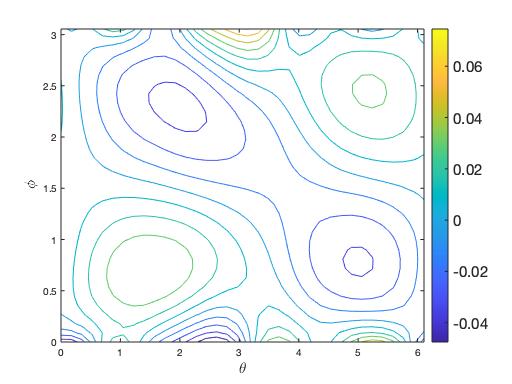}
		 			\label{fig:8k}
	\end{subfigure}%
	\vskip\baselineskip
	\caption{Heat inverse problem on a semi-torus:  $u_0(\alpha, \beta) = 2 + \sin(\alpha)\cos(\beta)$ with $\tau = 0.012$, $t^* = 5$. }
	 		\label{4}
\end{figure}

 % The primary reason of rather two distinct results for elliptic and heat inverse problem is that, in the latter case, one can control the degree of ill-posedness by reducing the time at which the solution is observed.

\section{Conclusions and Open Directions}\label{sec:conclusions}

In this paper, we developed graph-based Mat\'ern priors for solving Bayesian inverse problems on manifolds with boundaries. Our idea is to extend the Mat\'ern priors introduced in \cite{harlim2020kernel}, developed for elliptic PDEs on closed manifolds, by representing the boundary conditions via a set of functions obtained from solving Laplace equations on manifolds with appropriate Dirichlet boundary conditions. To solve PDEs on manifolds with boundaries, we employed the recently developed GPDM algorithm \cite{jiang2020ghost}, which uses fictitious ghost points to remove the bias induced by integrating radial type kernels near the boundaries. 

We validated this approach on two test problems. The first problem is an inversion of the diffusion coefficient of an elliptic PDE from the solution of the PDE corrupted by noise. The second problem is an inversion of the initial condition of a heat equation from noisy observation of the solution at a positive time. Based on our numerical simulations, we found positive results given the ill-posedness of the inverse problems we considered.
%with a credible interval when the problems are not severely ill-posed. In a severely ill-posed case, we showed that we can recover the PDE solution accurately from a series of diffusion coefficients whose spatial profiles are completely different than the underlying coefficients.  

While the proposed method produces encouraging results, there are many open questions. First, we should point out that while the computational cost is independent of the ambient dimension, it scales exponentially as a function of intrinsic dimension (see \cite{jiang2020ghost} for detailed convergence rates for the forward maps). Since the Bayesian inversion involves evaluations of the forward maps many times on the proposals drawn from the prior, it is desirable to have a computationally cheaper forward map. Second, the method represents the hidden variables by a vector whose components are the function values of the variable of interest (e.g., diffusion coefficients) on the given point clouds. How to extend this to other points on the domain is of practical interest. Beyond these practical considerations, it is also of interest to understand the theoretical aspect of such an approach, especially the effects of boundaries, extending the theoretical convergence result for closed manifolds in \cite{harlim2020kernel}.

\section*{Acknowledgment}
The research of JH was partially supported under the NSF grant DMS-1854299. DSA is thankful for the support of NSF and NGA through the grant DMS-2027056. The work of DSA was also partially supported by the NSF Grant DMS-1912818/1912802.

\bibliographystyle{abbrvnat} 
\bibliography{references}

\end{document}